\documentstyle[amscd]{amsart}
\numberwithin{equation}{section}
\theoremstyle{plain}
\newtheorem{thm}{Theorem}[section]
\newtheorem{cor}[thm]{Corollary}
\newtheorem{lem}[thm]{Lemma}
\newtheorem{prop}[thm]{Proposition}
\begin{document}
\title{Topological full groups of  $C^*$-algebras
arising from $\beta$-expansions}
\author{Kengo Matsumoto}
\address{ 
Department of Mathematics, 
Joetsu University of Education,
Joetsu 943-8512 Japan}
\author{Hiroki Matui}
\address{ 
Graduate School of Science, 
Chiba University,
Inage-ku, Chiba 263-8522 Japan}
\maketitle
\begin{abstract}
We will introduce a family $\Gamma_\beta, 1 < \beta \in {\mathbb{R}}$ 
of infinite non-amenable discrete groups as an interpolation of the
Higman-Thompson groups $V_n, 1 < n \in {\mathbb{N}}$
by using the topological full groups of the groupoids defined by 
$\beta$-expansions of real numbers.
They are regarded as full groups of certain interpolated Cuntz algebras.
The groups $\Gamma_\beta, 1 < \beta \in {\mathbb{R}}$ 
are realized as groups of piecewise linear functions on $[0,1]$
if the $\beta$-expansion of $1$ is finite or ultimately periodic.
We also classify the groups $\Gamma_\beta, 1 < \beta \in {\mathbb{R}}$
by the number theoretical property of $\beta$.    
\end{abstract}

\def\Zp{{ {\mathbb{Z}}_+ }}
\def\Ext{{{\operatorname{Ext}}}}
\def\Im{{{\operatorname{Im}}}}
\def\Hom{{{\operatorname{Hom}}}}
\def\Aut{{{\operatorname{Aut}}}}
\def\Ker{{{\operatorname{Ker}}}}
\def\dim{{{\operatorname{dim}}}}
\def\det{{{\operatorname{det}}}}
\def\supp{{{\operatorname{supp}}}}
\def\id{{{\operatorname{id}}}}
\def\OB{{ {\mathcal{O}}_\beta }}
\def\A{{ {\mathcal{A}}}}
\def\AB{{ {\mathcal{A}}_\beta }}
\def\DB{{ {\mathcal{D}}_\beta }}
\def\FB{{ {\mathcal{F}}_\beta }}
\def\E{{ {\mathcal{E}} }}



\section{Introduction}

The class of finitely  presented  non-abelian infinite groups 
is one of the most interesting and important classes of infinite groups from the viewpoints of not only group theory but also geometry and topology.
The study of finitely presented simple infinite groups has begun with Richard J. Thompson in 1960's. 
He has discovered first two such groups in \cite{Thompson}.
They are now known as the groups $V_2$ and $T_2$. 
G. Higman and K. S. Brown have generalized his examples to
infinite family of finitely presented infinite groups.
One of such family is the groups written  
$V_n, 1<n \in {\mathbb{N}}$ which are called the Higman-Thompson groups.
They are all finitely presented 
and their commutator subgroups
are all simple.
Their abelianizations are trivial if $n$ is even, and ${\mathbb{Z}}_2$ if $n$ is odd.
The Higman-Thompson group 
$V_n$ 
is represented as the group of right continuous piecewise linear functions 
$f: [0,1) \longrightarrow [0,1)$
having finitely many singularities such that
all singularities of $f$ are in ${\mathbb{Z}}[\frac{1}{n}]$,
the derivative of $f$ at any non-singular point is
$n^k$ for some $k \in {\mathbb{Z}}$ and 
$f({\mathbb{Z}}[\frac{1}{n}] \cap [0,1)) 
= {\mathbb{Z}}[\frac{1}{n}] \cap [0,1)$
(\cite{Thompson}).
 V. Nekrashevych \cite{Nek} has shown that 
the Higman-Thompson group $V_n$ appears as  a certain subgroup
of the unitary group of the Cuntz algebra ${\mathcal{O}}_n$ of order $n$.
The subgroup
of the unitary group of ${\mathcal{O}}_n$
is the 
continuous full group $\Gamma_n$
of ${\mathcal{O}}_n$,
which is also called the topological full group of the associated groupoid
 (see also \cite[Remark 6.3]{MatuiPre2012}).
Recently the authors have independently studied full groups of 
the Cuntz-Krieger algebras 
and full groups of the groupoids coming from shifts of finite type.
The first 
named author has studied the normalizer groups of the Cuntz-Krieger algebras
which are called the continuous full groups from the view point of orbit equivalences of topological Markov shifts and $C^*$-algebras 
(\cite{MaPacific}, \cite{MaPAMS}, etc.).
He has proved that the continuous full groups are complete invariants for 
the continuous orbit equivalence classes of the underlying topological Markov shifts.
The second named author has studied the continuous full groups of more general  
\'{e}tale groupoids (\cite{MatuiPLMS},
\cite{MatuiPre2011}, \cite{MatuiPre2012}, etc.).
He has called them the topological full groups of \'{e}tale groupoids.
He has proved that if 
an \'{e}tale groupoid is minimal, 
the topological full group of the groupoid is a complete invariant 
for the isomorphism  class of the groupoid.
He has also shown that if a groupoid comes from a shift of finite type,
the topological full group is of type $F_\infty$ and in particular finitely presented. 
He has furthermore obtained that 
the topological full groups for shifts of finite type 
are simple if and only if the homology group $H_0(G)$ of the groupoid $G$
is $2$-divisible.
Hence we know an infinite family of finitely presented infinite simple groups 
coming from symbolic dynamics. 
V.  Nekrashevych's paper \cite{Nek} says that 
the Higman-Thompson groups appear 
as the topological full groups 
of the groupoids of the full shifts and 
as the continuous full groups 
of the Cuntz algebras.
In \cite{KMW},
a family of $C^*$-algebas $\OB, 1 <\beta \in {\mathbb{R}}$
has been introduced. 
It arises from a family of certain subshifts called 
the $\beta$-shifts which are the symbolic dynamics 
defined by the $\beta$-transformations on the unit interval $[0,1]$.
The family of the $\beta$-shifts 
is an interpolation of the full shifts.
Hence the $C^*$-algebras    
$\OB, 1 <\beta \in {\mathbb{R}}$
are considered as interpolation of the Cuntz algebras
${\mathcal{O}}_N, 1 <N \in {\mathbb{N}}$.

In the present paper, 
we will introduce a family 
$\Gamma_\beta, 1 <\beta \in {\mathbb{R}}$
of infinite discrete groups 
as an interpolation of the Higman-Thompson groups 
$V_n, 1<n \in {\mathbb{N}}$
such that $\Gamma_n = V_n, 1<n \in {\mathbb{N}}.$
The  groups
 $\Gamma_\beta, 1<\beta \in {\mathbb{R}} $ 
are defined as the continuous full groups of the 
 $C^*$-algebras    
$\OB, 1 <\beta \in {\mathbb{R}}$.
They are also considered as the topological full groups 
of the 
\'{e}tale groupoids $G_\beta$ for the $\beta$-shifts.
We will first study the groupoid $G_\beta$ and show that 
the groupoid $G_\beta$ for each $1<\beta \in {\mathbb{R}}$
is an essentially principal, purely infinite, minimal \'{e}tale groupoid.
The homology groups $H_i(G_\beta)$ are computed as
\begin{equation}
H_i(G_\beta) \cong
\begin{cases}
K_i(\OB) & \text{ if } i=0,1, \\
0 & \text{ if } i\ge 2. 
\end{cases} 
\end{equation}

We will show that

\begin{thm}[Theorem \ref{thm:Theorem1}]
Let $1< \beta \in {\mathbb{R}}$ be a real number.
Then 
the group
$\Gamma_\beta$ is a countably infinite discrete non-amenable 
group such that
its commutater subgroup
$D(\Gamma_\beta)$ is simple.
\end{thm}
For a real number $\beta >1$, let us denote by 
$d(1,\beta)=\xi_1 \xi_2 \xi_3 \cdots $
the $\beta$-adic expansion of $1$ 
which means 
$\xi_i \in {\mathbb{Z}}$, 
$0 \le \xi_i \le [\beta]$ and
\begin{equation*}
1 = \frac{\xi_1}{\beta} + \frac{\xi_2}{\beta^2} +\frac{\xi_3}{\beta^3} 
    + \cdots. 
\end{equation*}
The expansion 
$d(1,\beta)$
is said to be finite if there exists $k\in {\mathbb{N}}$
such that 
$\xi_{m} =0$ for all $m >k$.
If there exists $l\le k$ such that 
$$
d(1,\beta)=\xi_1 \cdots \xi_l 
\xi_{l+1} \cdots \xi_{k+1}
\xi_{l+1} \cdots \xi_{k+1}
\xi_{l+1} \cdots \xi_{k+1}
\cdots,
$$
the expansin 
$d(1,\beta)$ is said to be ultimately periodic 
and written
$d(1,\beta) = $
$\xi_1 \cdots \xi_l$
$\Dot{\xi}_{l+1} \cdots \Dot{\xi}_{k+1}$.
It is well-known that
the Higman-Thompson group 
$V_n, n \in {\mathbb{N}}$ 
is represented as the group of right continuous piecewise linear functions 
$f: [0,1) \longrightarrow [0,1)$
having finitely many singularities such that
all singularities of $f$ are in ${\mathbb{Z}}[\frac{1}{n}]$,
the derivative of $f$ at any non-singular point is
$n^k$ for some $k \in {\mathbb{Z}}$ and 
$f({\mathbb{Z}}[\frac{1}{n}] \cap [0,1)) 
= {\mathbb{Z}}[\frac{1}{n}] \cap [0,1)$.
We will introduce a notion of $\beta$-adic PL functions on the interval $[0,1]$,
and show the following.
\begin{thm}[Theorem \ref{thm:SFTPL} and Theorem \ref{thm:soficPL}]
Let $1< \beta \in {\mathbb{R}}$ be a real number
such that the 
$\beta$-expansion $d(1,\beta)$ of $1$ is finite or ultimately periodic.
Then the group
$\Gamma_\beta$ is realized as the group of $\beta$-adic PL functions on the interval
$[0,1]$.
\end{thm}

It is well-known
that
$d(1,\beta)$
is  finite if and only if 
the $\beta$-shift $(X_\beta, \sigma)$ is a shift of finite type,
and 
$d(1,\beta)$
is  ultimately periodic if and only if 
the $\beta$-shift $(X_\beta, \sigma)$ is a sofic shift.
If $\beta = \frac{1 + \sqrt{5}}{2}$,
the number is the positive solution of the quadratic equation
$\beta^2 = \beta +1$ so that the $\beta$-expansion is finite
$d(1,\beta) = 110000\cdots$. 
We will classify the interpolated Higman-Thompson groups 
$\Gamma_\beta, 1 <\beta \in {\mathbb{R}}$
in the number theoretical property of $\beta$
and the dynamical property of the $\beta$-shift $(X_\beta,\sigma)$
in the following way.
\begin{thm}[Theorem \ref{thm:HTfinite}, Theorem \ref{thm:HTsofic} 
and Theorem \ref{thm:HTnonsofic}]
Let $1 <\beta \in {\mathbb{R}}$ be a real number
and
$d(1,\beta)= \xi_1 \xi_2 \xi_3 \cdots$ 
be the $\beta$ expansion of $1$.
\begin{enumerate}
\renewcommand{\labelenumi}{(\roman{enumi})}
\item
If
$d(1,\beta)$ is finite, that is, 
$d(1,\beta) = \xi_1 \xi_2 \cdots \xi_k 00\cdots$,
then
the group
$\Gamma_\beta$ is isomorphic to the Higman-Thompson group
$V_{\xi_1+ \cdots +\xi_k +1}$
of order   
$\xi_1+ \cdots +\xi_k +1$.
\item
If $d(1,\beta)$ is ultimately periodic, that is, 
$d(1,\beta) = \xi_1 \cdots \xi_l \Dot{\xi}_{l+1} \cdots \Dot{\xi}_{k+1}$,
then 
the group
$\Gamma_\beta$ is isomorphic to the Higman-Thompson group
$V_{\xi_{l+1}+ \cdots +\xi_{k +1}}$
of order
$\xi_{l+1} +\cdots +\xi_{k +1}$.
\item
If $1 <\beta \in {\mathbb{R}}$
is not ultimately periodic,
then the group
$\Gamma_\beta$ is not isomorphic to any of 
the Higman-Thompson group $V_n, n \in {\mathbb{N}}$. 
\end{enumerate}
\end{thm}

\section{Preliminaries of the $C^*$-algebra $\OB$}

We fix an arbitrary real number $\beta > 1$ unless we specify.
Take a natural number $N$ with $N-1 < \beta \le N$.
Put $\Sigma = \{ 0,1,...,N-1 \}$. 
For a non-negative real number $t$,
 we denote by $[t]$ the integer part of $t$. 
Let $f_\beta : [0,1] \rightarrow [0,1]$ be the funtion defined by 
\begin{equation*}
f_\beta (x) = \beta x - [\beta x], \qquad  x \in [0,1].
\end{equation*} 
The $\beta$-expansion of $x \in [0,1]$ is a sequence 
$\{ d_n(x,\beta)\}_{n \in {\mathbb N}}$ 
of integers of $\Sigma$ determined by  
\begin{equation*}
       d_n(x,\beta) = [\beta f_{\beta}^{n-1}(x)], \qquad n \in {\mathbb{N}}
\qquad (cf.  \cite{Parry}, \cite{Renyi}).
\end{equation*}
The numbers 
$ d_n(x,\beta) $ will be denoted by $d_n(x)$ for simplicity.
By this sequence, we can write $x$ as
\begin{equation*}
        x = \sum_{n=1}^{\infty} \frac{d_n(x)}{\beta^n}.
\end{equation*}

We endow the infinite product $\Sigma^{\mathbb{N}}$ 
with the product topology and the lexicographical order
from the leftmost letter 
$x_1$ in 
$(x_1,x_2,\dots) = (x_n)_{n \in {\mathbb{N}}} \in \Sigma^{\mathbb{N}}$. 
We denote by $\sigma$ the shift on $\Sigma^{\mathbb{N}}$ defined by 
$
\sigma((x_n)_{n \in {\mathbb{N}}}) = (x_{n+1})_{n \in {\mathbb{N}}}.
$
Let $\xi_\beta =(\xi_n)_{n \in {\mathbb N}} \in \Sigma^{\mathbb{N}}$ 
be the supremum element of 
$\{ (d_n(x))_{n \in {\mathbb{N}}} \mid x \in [0,1)\}$ 
with respect to the lexicographical order in $\Sigma^{\mathbb{N}}$
which is defined by
\begin{equation*}
\xi_\beta = \sup_{x \in [0,1)} (d_n(x))_{n \in {\mathbb N}}.
\end{equation*}
Define the $\sigma$-invariant compact subset $X_{\beta}$ 
of $\Sigma^{\mathbb{N}}$ by 
\begin{equation*}
  X_{\beta} = \{ \omega \in \Sigma^{\mathbb{N}} | 
               \sigma^m(\omega) \le \xi_{\beta}, m= 0, 1, 2,... \}.
\end{equation*} 

\noindent
{\bf Definition (cf. \cite{Parry}, \cite{Renyi}).} 
The subshift 
$(X_{\beta}, \sigma)$ 
is called the $\beta$-shift.
 
\medskip

Example 1. $\beta = N \in {\mathbb{N}}$ with $N>1$.
 As 
$\xi_{\beta} = (N-1)(N-1)\cdots $,
the subshift
$$
  X_{N} = \{ (x_n)_{n \in {\mathbb N}} \in \prod_{n=1}^{\infty}\{0,1,\dots, N-1\} \mid
              x_n= 0,1,\dots, N-1\}
$$
is the full $N$-shift.

Example 2. $\beta = \frac{1+ \sqrt{5}}{2}$. 
As $N=2$ and $d(1,\beta) = 1100\cdots, \xi_{\beta}= 10101010 \cdots$,
we see 
\begin{equation*}
  X_{\frac{1+ \sqrt{5}}{2}} = \{ (x_n)_{n \in {\mathbb N}} \in \prod_{n=1}^{\infty}\{0,1\} \mid
    \text{\lq\lq} 11 \text{\rq\rq} 
    \text{  does not appear in } (x_n)_{n \in {\mathbb N}} \}.
\end{equation*}
This  is a shift of finite type 
 $X_A$ determined by the matrix 
\begin{math}
A=
\bigl[
\begin{smallmatrix}
1 & 1 \\
1 & 0
\end{smallmatrix}
\bigr]
\end{math}.

Example 3. $\beta = 2+ \sqrt{3}$. 
As $N=4$ and $d(1,\beta) = \xi_{\beta}= 3 \Dot{2}$,
we see 
\begin{equation*}
  X_{2+ \sqrt{3}} = \{ (x_n)_{n \in {\mathbb N}} \in \prod_{n=1}^{\infty}\{0,1,2,3\} \mid
     (x_{n+m})_{n \in {\mathbb N}} \le 3\Dot{2} \text{ for all } m=0,1,2,\dots \}.
\end{equation*}
This  is a sofic shift but not a shift of finite type.

Example 4. $\beta = \frac{3}{2}$. 
As $N=2$ and $\xi_{\beta}= 101000001 \cdots$, 
we see 
$$
  X_{\frac{3}{2}} = \{ (x_n)_{n \in {\mathbb N}} \in \prod_{n=1}^{\infty}\{0,1\} |
     (x_{n+m})_{n \in {\mathbb N}} \le 101000001 \cdots, m=0,1,2, \dots
 \}.
$$
This is not a sofic shift
 (hence not a shift of finite type).

A finite sequence 
$\mu = (\mu_1, \dots, \mu_k) $  
of elements $\mu_j \in \Sigma$
is called 
 a block or a word. 
We denote by $|\mu|$ the length $k$ of $\mu$.
Set for $k \in {\mathbb{N}}$
\begin{equation*}
 B_k(X_{\beta}) = \{ \mu | \text{ a block with length } k 
                \text{ appearing in some } x \in X_{\beta} \} 
\end{equation*}
and
$ 
B_*(X_{\beta})  = \bigcup_{k=0}^\infty B_k(X_{\beta})
$ 
where  
$B_0(X_{\beta})$ 
denotes the empty word $\emptyset$.

In \cite{KMW},
a family $\OB, 1< \beta \in {\mathbb{R}}$
of simple purely infinite $C^*$-algebras has been introduced 
as the $C^*$-algebras associated with $\beta$-shifts
$(X_\beta,\sigma)$.
For $\beta = N \in {\mathbb{N}}$,
the $C^*$-algebra is isomorphic to the Cuntz algebra ${\mathcal{O}}_N$.
Hence the family $\OB, 1< \beta \in {\mathbb{R}}$ 
is an interpolation of the Cuntz algebras
${\mathcal{O}}_N, 1< N \in {\mathbb{N}}$. 
As in \cite{KMW},
let us denote by
$S_0, S_1,\dots,S_{N-1}$
the generating partial isometries 
of $\OB$
which satisfy the equalities
\begin{equation*}
\sum_{j=0}^{N-1}S_j S_j^* =1, 
\qquad
S_i^*S_i = 1 \quad \text{ for } i=0,1,\dots,N-2.
\end{equation*}
We put
$S_\mu = S_{\mu_1}\cdots S_{\mu_l}$
for $\mu = (\mu_1,\dots,\mu_l) \in B_l(X_\beta)$
and
$
a_\mu = S_\mu^* S_\mu.
$
We define $C^*$-subalgebras
of $\OB$: 
\begin{align*}
{\mathcal{A}}_l
& = \text{the } C^*\text{-subalgebra of } \OB
\text{ generated by } 
S_\mu^*S_\mu, \mu \in B_l(X_\beta),\\
\AB
& = \text{the } C^*\text{-subalgebra of } \OB
\text{ generated by } 
S_\mu^*S_\mu, \mu \in B_*(X_\beta),\\
\DB
& = \text{the } C^*\text{-subalgebra of } \OB
\text{ generated by } 
S_\mu a S_\mu^*, \mu \in B_*(X_\beta), a \in \AB, \\
\FB
& = \text{the } C^*\text{-subalgebra of } \OB
\text{ generated by } 
S_\mu a S_\nu^*, \mu,\nu  \in B_k(X_\beta),
k \in \Zp, a \in \AB.
\end{align*}
As
$S_\mu^*S_\mu = S_\mu^*S_0^*S_0S_\mu$,
the algebra 
${\mathcal{A}}_l$ 
is naturally embedded into 
${\mathcal{A}}_{l+1}$.
It is commutative and finite dimensional
so that
the algebras $\AB$, $\DB$ and $\FB$ are   
all AF-algebras, in particular
$\AB$ and  $\DB$ are both commutative.
Put
$\rho_j(x) = S_j^* x S_j$
for $x \in \AB, j=0,1,\dots,N-1$.
Then the $C^*$-algebra $\OB$ has a universal property 
subject to the relations:
\begin{equation*}
\sum_{j=0}^{N-1}S_j S_j^* =1, 
\qquad
\rho_j(x) = S_j^* x S_j
 \quad \text{ for } 
x \in \AB, \, j=0,1,\dots,N-1 \qquad (cf. \cite{MaDocMath02}).
\end{equation*}
For $t \in {\mathbb{R}}/{\mathbb{Z}} ={\mathbb{T}}$
the correspondence
$S_j \longrightarrow e^{2\pi \sqrt{-1} t} S_j, j=0,1,\dots,N-1$
yields an automorphism of $\OB$
which 
gives rise to an action on $\OB$ of ${\mathbb{T}}$
called the gauge action written $\hat{\rho}$.
The gauge action has a unique KMS-state 
denoted by $\varphi$
on $\OB$ at 
inverse temperature $\log \beta$.
For
the projections
$
a_{\xi_1 \cdots \xi_n} = 
S_{\xi_1 \cdots \xi_n}^*S_{\xi_1 \cdots \xi_n} \in \AB, 
n=1,2, \dots,
$ 
the values 
$\varphi(a_{\xi_1 \cdots \xi_n})$
are computed as 
\begin{equation*}
\varphi(a_{\xi_1 \cdots \xi_n}) = 
\beta^n - \xi_1 \beta^{n-1} - \cdots -\xi_{n-1}\beta -\xi_n 
= \sum_{i=1}^\infty \frac{\xi_{i+n}}{\beta^i}, \qquad
n=1,2,\dots (\cite{KMW}).
\end{equation*}
Let $m(l)$ denote the dimension $\dim {\mathcal{A}}_l$
of ${\mathcal{A}}_l$.
Denote by
$E_1^l,\dots,E_{m(l)}^l$ 
the set of minimal projections of 
 ${\mathcal{A}}_l$.
As in \cite[Lemma 3.3]{KMW},
the projection $E_i^l, i=1,\dots,m(l)$
is of the form
$E_i^l = a_{\xi_1\cdots\xi_{p_i}}
-
a_{\xi_1\cdots\xi_{q_i}}
$
for some 
$p_i, q_i =0,1,\dots$.
The projections
$a_{\xi_1\cdots\xi_n}, n \in \Zp$ 
are totally ordered by the value 
$\varphi(a_{\xi_1\cdots\xi_n})$.
We order 
$E_1^l, \dots,E_{m(l)}^l$ 
following the order
$
   \varphi(a_{\xi_1\cdots\xi_{p_1}}) 
< \cdots
< \varphi(a_{\xi_1\cdots\xi_{p_{m(l)}}})$
in ${\mathbb{R}}$.

Some basic subclasses of $\beta$-shifts are characterized 
in terms of the $\beta$-expansion $d(1,\beta)$ of $1$ and the projections
$a_{\xi_1 \cdots \xi_n}$
in the following way. 
\begin{lem}[{\cite{Parry}, cf. \cite[Proposition 3.8]{KMW}}]
The followings are equivalent:
\begin{enumerate}
\renewcommand{\labelenumi}{(\roman{enumi})}
\item
$(X_\beta,\sigma)$ is a shift of finite type.
\item
$d(1,\beta)$ is finite, that is, 
$d(1,\beta) = \xi_1 \xi_2 \cdots \xi_k 000\cdots$
for some $k \in {\mathbb{N}}$. 
\item
 $a_{\xi_1 \cdots \xi_k} =1$ for some $k \in {\mathbb{N}}.$ 
\end{enumerate}
\end{lem}
We call $(X_\beta,\sigma)$ an SFT $\beta$-shift if 
$(X_\beta,\sigma)$ is a shift of finite type.
\begin{lem}[{\cite{BM}, cf. \cite[Proposition 3.8]{KMW}}]
The followings are equivalent:
\begin{enumerate}
\renewcommand{\labelenumi}{(\roman{enumi})}
\item
$(X_\beta,\sigma)$ is a sofic shift. 
\item
$d(1,\beta)$ is ultimately periodic, that is, 
$d(1,\beta) = \xi_1 \cdots \xi_l \Dot{\xi}_{l+1} \cdots \Dot{\xi}_{k+1}$
for some $l\le k$.
\item
 $a_{\xi_1 \cdots \xi_l} =a_{\xi_1 \cdots \xi_{k+1}}$
 for some $l \le k \in {\mathbb{N}}.$
\end{enumerate}
\end{lem}
We call $(X_\beta,\sigma)$ a sofic $\beta$-shift if 
$(X_\beta,\sigma)$ is a sofic shift.

The K-groups of the $C^*$-algebra $\OB$ have been computed 
in the following way. 
\begin{lem}[{\cite{KMW}}]
\begin{align*}
K_0(\OB) 
=
& 
{
\begin{cases}
{\mathbb{Z}}/(\xi_1+ \xi_2+ \cdots + \xi_k -1){\mathbb{Z}}
& \text{ if } d(1,\beta) = \xi_1 \xi_2 \cdots \xi_k 000\cdots, \\
{\mathbb{Z}}/(\xi_{l+1}+ \cdots +\xi_{k+1}){\mathbb{Z}}
& \text{ if } d(1,\beta) = \xi_1 \cdots \xi_l 
\Dot{\xi}_{l+1} \cdots \Dot{\xi}_{k+1}, \\
{\mathbb{Z}} & \text{ otherwise.  } 
\end{cases}
} \\
\intertext{The position  $[1]$  of the unit of $\OB$
in $K_0(\OB)$ corresponds to the class 
$[1]$ of $ 1 \in {\mathbb{Z}} $
in the first two cases, and to
$1 \in {\mathbb{Z}}$
 in the third case, and}\\ 
 K_1(\OB)
= & 0 \qquad \text{ for any } \beta >1.
\end{align*}
\end{lem}

\section{Topological full groups of the groupoid $G_\beta$} 
The $C^*$-algebra 
$\OB, 1< \beta \in {\mathbb{R}}$
has been originally constructed as 
the $C^*$-algebra associated with subshift
$(X_\beta,\sigma), 1< \beta \in {\mathbb{R}}$.
It is regarded as the  $C^*$-algebra 
$C^*_r(G_\beta)$
of 
a certain essentially principal \'{e}tale groupoid 
$G_\beta$ as in \cite{MaDocMath02}.
Denote by 
$G_\beta^{(0)}$ the unit space of the groupoid
$G_\beta$, so that 
the $C^*$-algebra
$C(G_\beta^{(0)})$
of all complex valued continuous functions on 
$G_\beta^{(0)}$ is canonically isomorphic to 
the $C^*$-subalgebra $\DB$ of $\OB$
which is called the canonical Cartan subalgebra of $\OB$.
There exists a continuous surjection
$\sigma_\beta$ on $G_\beta^{(0)}$
such that
\begin{equation*}
G_\beta = 
\{(x,k-l,y) \in G_\beta^{(0)} \times {\mathbb{Z}} \times G_\beta^{(0)} 
\mid
\sigma_\beta^k(x) = \sigma_\beta^l(y)
\text{ for some } k,l \in \Zp\}.
\end{equation*}

For an \'{e}tale groupoid $G$ 
we let  
$G^{(0)}$ denote the unit space of $G$ and let
$s$ and $r$ denote the source map and range map.
For $x \in G^{(0)}$, the set 
$G(x) =r(Gx)$ is called the $G$-orbit of $x$.
If every $G$-orbit is dense in $G^{(0)}$, 
$G$ is said to be minimal (\cite{MatuiPre2012}, \cite{Renault}). 
\begin{lem}
For $1<\beta \in {\mathbb{R}}$,
the groupoid $G_\beta$ is an essentially principal, minimal groupoid.
\end{lem}
\begin{pf}
The $C^*$-subalgebra $\FB$ of $\OB$ is the $C^*$-algebra
$C^*_r(H_\beta)$ of an AF-subgroupoid $H_\beta$ of $G_\beta$,
which is defined by 
\begin{equation*}
H_\beta = 
\{(x,0,y) \in G_\beta^{(0)} \times {\mathbb{Z}} \times G_\beta^{(0)} 
\mid
\sigma_\beta^k(x) = \sigma_\beta^k(y) \text{ for some } k \in \Zp \}.
\end{equation*}
 As the algebra $\FB$ is simple (\cite[Proposition 3.5]{KMW}), 
the groupoid  $H_\beta$ 
is minimal so that $G_\beta$ is minimal.
\end{pf}
A subset $U \subset G$ is called a $G$-set if 
$r|_U, s|_U$ are injective.
The homeomorphism
$r \circ (s|_{U})^{-1}$ from $s(U)$ to $r(U)$
is denoted by $\pi_U$. 
Following \cite{MatuiPre2012},
$G$ is said to be purely infinite
if for every clopen set $A \subset G^{(0)}$
there exist clopen  $G$-sets $U,V \subset G$ such that 
$s(U) = s(V) = A, r(U) \cup r(V) \subset A, r(U) \cap r(V) =\emptyset.$  
\begin{lem}
For $1<\beta \in {\mathbb{R}}$,
the groupoid $G_\beta$ is purely infinite.
\end{lem}
\begin{pf}
As the $C^*$-algebra 
$\DB$ is isomorphic to
the algebra $C(G_\beta^{(0)})$
of continuous functions on 
$G_\beta^{(0)}$,
we may identify the projections of 
$\DB$ with the
clopen sets of $G_\beta^{(0)}$.
Hence a 
clopen set of $G_\beta^{(0)}$
may be considered as a finite sum of the form 
$P= S_\mu E_i^l S_\mu^*$
for some $\mu \in B_k(X_\beta)$ with $k \le l$ 
such that
$S_\mu^* S_\mu \ge E_i^l$.
It is enough to consider $P= S_\mu E_i^l S_\mu^*$ for simplicity.
The minimal projection
$E_i^l \in \A_l$ is of the form
$E_i^l = a_{\xi_1\cdots\xi_{p_i}}
-
a_{\xi_1\cdots\xi_{q_i}}
$
for some 
$1 \le p_i, q_i \le l$
with
$ a_{\xi_1\cdots\xi_{p_i}}
>
a_{\xi_1\cdots\xi_{q_i}}.
$
Note that
\begin{equation}
S_\mu^* S_\mu \ge a_{\xi_1\cdots\xi_{p_i}}.
\label{eqn:smaxi}
\end{equation}
There exists $\gamma=(\gamma_1,\dots,\gamma_r) \in B_*(X_\beta)$
such that
\begin{equation*}
 ({\xi_1,\dots,\xi_{p_i}}, \gamma_1,\dots,\gamma_r) \in B_*(X_\beta),
\qquad
 ({\xi_1,\dots,\xi_{q_i}}, \gamma_1,\dots,\gamma_r) \not\in B_*(X_\beta).
\label{eqn:gammaxi}
\end{equation*}
Put the words
\begin{equation*}
\zeta_1(m) =(\overbrace{0,\dots, 0}^{m}),
\qquad
\zeta_2(m) =(\overbrace{0, \dots, 0,1}^{m}).
\end{equation*}
By \cite[Corollary 3.2]{KMW},
there exists $m \in {\mathbb{N}}$ such that 
\begin{equation}
 a_{\xi_1\cdots\xi_{p_i} \gamma \zeta_1(m)}
=
 a_{\xi_1\cdots\xi_{p_i} \gamma \zeta_2(m)}
=1.
\label{eqn:axigammazeta}
\end{equation}
Put
$
\zeta_1 = \zeta_1(m),
\zeta_2 = \zeta_2(m).
$
By \eqref{eqn:smaxi}, \eqref{eqn:axigammazeta}, we have
\begin{equation*}
a_{\mu \gamma \zeta_1} 
\ge S_{\gamma \zeta_1}^* a_{\xi_1\cdots\xi_{p_i}} S_{\gamma \zeta_1} 
= a_{\xi_1\cdots\xi_{p_i} \gamma \zeta_1} =1   
\end{equation*}
so that
$a_{\mu \gamma \zeta_1} =1$ and similarly 
$a_{\mu \gamma \zeta_2} =1.$ 
We set 
\begin{equation*}
U = S_{\mu \gamma \zeta_1} E_i^l S_{\mu}^*,
\qquad
V = S_{\mu \gamma \zeta_2} E_i^l S_{\mu}^*,
\end{equation*}
which correspond to ceratin clopen $G$-sets in $G_\beta$.
It then follows that
\begin{equation*}
U^* U = S_{\mu} E_i^l a_{\mu \gamma \zeta_1} E_i^l S_{\mu}^*
= S_{\mu} E_i^l S_{\mu}^* =P
\quad
\text{ and similarly }
\quad
V^* V  =P
\end{equation*}
so that   
\begin{equation*}
U U^* + V V^*  
= S_{\mu \gamma \zeta_1} E_i^l S_{\mu \gamma \zeta_1}^*
+ S_{\mu \gamma \zeta_2} E_i^l S_{\mu \gamma \zeta_2}^*.
\end{equation*}
As
\begin{equation*}
S_{\gamma \zeta_1}^*
 E_i^l 
S_{\gamma \zeta_1}
=
S_{\zeta_1}^*S_{\gamma}^*
( a_{\xi_1\cdots\xi_{p_i}} - a_{\xi_1\cdots\xi_{q_i}})
S_{\gamma}S_{\zeta_1}
=
S_{\zeta_1}^* a_{\xi_1\cdots\xi_{p_i}\gamma}S_{\zeta_1}
 =1,
\end{equation*}
we have
\begin{equation*}
P S_{\mu \gamma \zeta_1} E_i^l S_{\mu \gamma \zeta_1}^*
=S_{\mu} E_i^l 
S_{\gamma \zeta_1} E_i^l S_{\mu \gamma \zeta_1}^*
=S_{\mu\gamma \zeta_1}
S_{\gamma \zeta_1}^*
 E_i^l 
S_{\gamma \zeta_1} E_i^l S_{\mu \gamma \zeta_1}^*
\end{equation*}
so that
$
P S_{\mu \gamma \zeta_1} E_i^l S_{\mu \gamma \zeta_1}^*
=S_{\mu\gamma \zeta_1}
 E_i^l S_{\mu \gamma \zeta_1}^*. 
$
This implies 
$UU^* \le P$ 
and similarly
$VV^*\le P$.
Since
$ S_{\mu\gamma \zeta_1} E_i^l S_{\mu\gamma \zeta_1}^*\cdot
  S_{\mu\gamma \zeta_2} E_i^l S_{\mu \gamma \zeta_2}^* =0,
$ 
we have
$UU^* + VV^* \le P$.
\end{pf}
Therefore we have
\begin{prop}\label{prop:pureinf}
For $1<\beta \in {\mathbb{R}}$,
the groupoid $G_\beta$ is an essentially principal, purely infinite, minimal \'{e}tale groupoid.
\end{prop}

We will next compute the homology groups $H_i(G_\beta)$
for the \'{e}tale groupoid $G_\beta$.
The homology theory for \'{e}tale groupoids has been studied in 
\cite{CraMoe}.
In \cite{MatuiPLMS}, the homology groups $H_i$
for the groupoids coming from shifts of finite type 
have been computed
such that the groups $H_i$ are isomorphic to the K-groups $K_i$ of the associated 
Cuntz-Krieger algebra for $i=0,1$, and $H_i =0$ for $i\ge 2$.
By following the argument of the proof of 
\cite[Theorem 4.14]{MatuiPLMS}, we have
\begin{prop}
For each $1 < \beta \in {\mathbb{R}}$,
the homology groups $H_i(G_\beta)$ are computed as
\begin{equation}
H_i(G_\beta) \cong
\begin{cases}
K_i(\OB) & \text{ if } i=0,1, \\
0 & \text{ if } i\ge 2. 
\end{cases} \label{eqn:Hi}
\end{equation}
\end{prop}
\begin{pf}
For each $1 < \beta \in {\mathbb{R}}$,
the map
$\rho_\beta: (x,n,y) \in G_\beta\longrightarrow n \in {\mathbb{Z}}$
gives rise to a groupoid homomorphism
such that the skew product 
$G_\beta \times_{\rho_\beta}{\mathbb{Z}}$ is  
homologically similar to the AF-groupoid $H_\beta$ 
(cf. \cite[Lemma 4.13]{MatuiPLMS}). 
We know that the groupoid $C^*$-algebra
$C^*_r(G_\beta \times_{\rho_\beta}{\mathbb{Z}})$
is stably isomorphic to
the crossed product $\OB\times_{\hat{\rho}}\mathbb{T}$
of $\OB$ by the gauge action, which is stably isomorphic to the AF-algebra
$C^*_r(H_\beta)$. 
Since the ${\mathbb{Z}}$-module structure on
 $H_0(G_\beta \times_{\rho_\beta}{\mathbb{Z}})$
is given by the induced action 
$\hat{\hat{\rho}}_*$
on $K_0(\OB\times_{\hat{\rho}}\mathbb{T})$
of the bidual action $\hat{\hat{\rho}}$
on $\OB\times_{\hat{\rho}}\mathbb{T}$,  
we get \eqref{eqn:Hi} 
by the same argument as \cite[Theorem 4.14]{MatuiPLMS}.
\end{pf}

In \cite{MatuiPLMS},
the notion of topological full groups for \'{e}tale groupoids has been introduced.
We will study the topological full groups of the groupoid
$G_\beta$ for the $\beta$-shift $(X_\beta,\sigma)$.

\noindent
{\bf Definition (\cite[Definition 2.3]{MatuiPLMS}).}
The topological full group
$[[G_\beta]]$ of the groupoid $G_\beta$ 
is defined by the group of all homeomorphisms $\alpha$ of $G_\beta^{(0)}$
such that $\alpha = \pi_U$ for some compact open $G_\beta$-set $U$. 
In what follows we denote the topological full group
$[[G_\beta]]$ by $\Gamma_\beta$.
By \cite[Proposition 5.6]{MatuiPLMS},
there exists a short exact sequence:
\begin{equation*}
1 \longrightarrow  
U(C(G^{(0)}_\beta) \longrightarrow 
N(C(G^{(0)}_\beta), C^*_r(G_\beta)) \longrightarrow 
\Gamma_\beta \longrightarrow  1
\end{equation*}
 where 
$U(C(G^{(0)}_\beta)$ denotes the group of unitaries in $C(G^{(0)}_\beta)$ 
and
$N(C(G^{(0)}_\beta), C^*_r(G_\beta)) $
denotes the group of unitaries in $C^*_r(G_\beta) $ which normalize 
$C(G^{(0)}_\beta)$.

Consider the full $n$-shift
$(X_n, \sigma)$ and its groupoid
$G_n$ (cf. \cite{MatuiPre2012}, \cite{Renault}). 
The groupoid $C^*$-algebra $C^*_r(G_n)$ is isomorphic to the Cuntz algebra 
${\mathcal{O}}_n$ of order $n$.
  V. Nekrashevych \cite{Nek} has shown that 
the Higman-Thompson group $V_n$ is identified with a certain subgroup
of the unitary group of ${\mathcal{O}}_n$.
The identification gives rise to an isomorphism
between the Higman-Thompson group $V_n$ and the topological full group
$\Gamma_n$ (see also \cite[Remark 6.3]{MatuiPre2012}).
 Hence our groups
 $\Gamma_\beta, 1<\beta \in {\mathbb{R}} $ 
are considered as interpolation of the Higman-Thompson group 
$V_n, 1<n \in {\mathbb{N}}$.
It is well-known that the groups
$V_n,  1< n \in {\mathbb{N}}$
are non-amenable and its commutator subgroups
$D(V_n)$ are all simple.
Proposition \ref{prop:pureinf} says that
the groupoid $G_\beta$ is an essentially principal, purely infinite, minimal groupoid
for every $1<\beta \in {\mathbb{R}}$.
By \cite[Proposition 4.10]{MatuiPre2012} and
\cite[Theorem 4.16]{MatuiPre2012},
we have a following generalization of the above fact for 
$V_n,   1< n \in {\mathbb{N}}$. 
\begin{thm}\label{thm:Theorem1}
Let $1 <\beta \in {\mathbb{R}}$ be a real number.
Then 
the group
$\Gamma_\beta$ is countably infinite discrete non-amenable 
group such that
its commutater subgroup
$D(\Gamma_\beta)$ is simple.
\end{thm}

\section{Realization of $\OB$ on $L^2([0,1])$}
The Higman-Thompson group 
$V_n, 1 < n \in {\mathbb{N}}$ 
is represented as the group of right continuous piecewise linear bijective functions 
$f: [0,1) \longrightarrow [0,1)$
having finitely many singularities such that
all singularities of $f$ are in ${\mathbb{Z}}[\frac{1}{n}]$,
the derivative of $f$ at any non-singular point is
$n^k$ for some $k \in {\mathbb{Z}}$ and 
$f({\mathbb{Z}}[\frac{1}{n}] \cap [0,1)) 
= {\mathbb{Z}}[\frac{1}{n}] \cap [0,1)$.
In order to represent our group
$\Gamma_\beta$ 
as a group of piecewise linear functions on $[0,1)$,
we will represent the algebra $\OB$ on $L^2([0,1])$
in the following way.

We denote by $H$
  the Hilbert space 
$L^2([0,1])$ of the square integrable functions on 
$[0,1]$ with respect to the Lebesgue measure.  
The essentially bounded measurable functions $L^\infty([0,1])$
 act on $H$ by multiplication.
We put the sequence
\begin{equation*}
\beta_n 
= 
\beta^n - \xi_1 \beta^{n-1} - \cdots -\xi_{n-1}\beta -\xi_n 
= \sum_{i=1}^\infty \frac{\xi_{i+n}}{\beta^i}, \qquad
n=1,2,\dots.
\end{equation*}
Consider the functions $g_0, g_1, \dots, g_{N-1}$ defined by 
\begin{align}
g_i(x)        &  =\frac{1}{\beta}(x + i)  \qquad  \text{for } i=0,1,\dots, N-2, \quad x \in [0,1], \\
g_{N-1}(x)  & =\frac{1}{\beta}(x + N-1) \qquad  \text{for } x \in [0,\beta_1]. 
\end{align}
They satisfy the following:
\begin{align}
\cup_{i=0}^{N-2}g_i([0,1]) & \cup g_{N-1}([0,\beta_1]) = [0,1], \\
f_\beta (g_i (x)) = x & \qquad \text{for } i=0,1,\dots, N-2, \quad x \in [0,1], \\
f_\beta(g_{N-1}(x)) = x & \qquad \text{for }  x \in [0,\beta_1].
\end{align}
For a measurable subset $U$ of $[0,1]$,
denote by 
$\chi_U$ the multiplication operator on $H$ 
of the characteristic function of $U$. 
Define  the bounded linear operators
$T_{f_\beta}$, $T_{g_i},i= 0,1,\dots,N-2 $
on $H$ by
\begin{align*}
(T_{f_\beta} \xi)(x) 
&  = \xi(f_\beta(x)) \qquad \text{for } \xi \in H, x \in [0,1], \\
(T_{g_i} \xi)(x) 
& = \xi(g_i(x))  \qquad \text{for }
 \xi \in H, x \in [0,1], \, i=0,1,\dots,N-2. 
\end{align*}
For the function $g_{N-1}$ on $[0,\beta_1]$
define the operator
$T_{g_{N-1}}$ by 
\begin{equation*}
(T_{g_{N-1}}\xi)(x) =
\begin{cases}
\xi( g_{N-1}(x)) & \text{ for  } x \in [0,\beta_1],\\
0 &  \text{ for  }  x\in (\beta_1,1].
\end{cases}
\end{equation*}
The followings are straightforward:
\begin{lem} Keep the above notations. 
We have
\begin{enumerate}
\renewcommand{\labelenumi}{(\roman{enumi})}
\item 
$T_{f_\beta}^* = \frac{1}{\beta}\sum_{i=0}^{N-1} T_{g_i}.$
\item 
$T_{f_\beta}^* T_{f_\beta} 
= \frac{N-1}{\beta} + \frac{1}{\beta} \chi_{[0,\beta_1]}.$
\item 
$T_{g_i}^* T_{g_i} =
\begin{cases}
\beta \chi_{[\frac{i}{\beta}, \frac{i+1}{\beta})} & \text{ for }  i=0,1,\dots,N-2,\\
\beta \chi_{[\frac{N-1}{\beta}, 1)} & \text{ for } i=N-1.
\end{cases}
$
\item 
$T_{g_i} T_{g_i}^* = 
\begin{cases}
\beta 1 & \text{ for }  i=0,1,\dots,N-2,\\
\beta \chi_{[0,\beta_1]} & \text{ for } i=N-1.
\end{cases}
$
\end{enumerate}
\end{lem}
We define
the operators $s_i, i=0,\dots,N-1$ on $H$ by setting
 \begin{equation*}
 s_i = \frac{1}{\sqrt{\beta}} T_{g_i}^*, \qquad i=0,1,\dots,N-1.
 \end{equation*}
By the above lemma, we have
\begin{prop}
The operators
$s_i, i=0,\dots,N-1$
are  partial isometries such that
\begin{align*}
s_i^* s_i & =
{\begin{cases}
1 & \text{ for }  i=0,1,\dots,N-2,\\
\chi_{[0,\beta_1]} & \text{ for } i=N-1,
\end{cases}} \\
s_i s_i^* & =    
{\begin{cases}
\chi_{[\frac{i}{\beta}, \frac{i+1}{\beta})} & \text{ for }  i=0,1,\dots,N-2,\\
\chi_{[\frac{N-1}{\beta}, 1)} & \text{ for } i=N-1
\end{cases}} \label{eqn:sisi} 
\quad
\text{ and hence }
 \quad
\sum_{i=0}^{N-1} s_is_i^* = 1.
\end{align*}      
\end{prop}
The natural ordering of
$\Sigma = \{ 0,1,\dots, N-1\}$
induces the lexicographical order
on $B_*(X_\beta)$, 
which means that 
for 
$\mu = (\mu_1,\dots,\mu_n) \in B_n(X_\beta)$
and
$\nu = (\nu_1,\dots,\nu_m) \in B_m(X_\beta)$,
the order 
$\mu \prec \nu$
is defined  
if $\mu_1 < \nu_1$ 
or $\mu_i = \nu_i$ for $i=1,\dots,k-1$ for some $k \le m,n$
and
$\mu_k < \nu_k$.
For a word
$\mu = (\mu_1,\dots,\mu_n) \in B_n(X_\beta)$,
denote by 
$\tilde{\mu} = (\tilde{\mu}_1,\dots,\tilde{\mu}_n) \in B_n(X_\beta)$,
the least word in $B_n(X_\beta)$ 
satisfying
$ (\mu_1,\dots,\mu_n) \prec (\tilde{\mu}_1,\dots,\tilde{\mu}_n)$.
If 
 $\mu = (\mu_1,\dots,\mu_n)$ is maximal in $B_n(X_\beta)$,
we set
$\tilde{\mu} = \emptyset.$
We will use the following notations for  
$\mu = (\mu_1,\dots,\mu_n) \in B_n(X_\beta)$,
\begin{equation*}
l(\mu) : = \frac{\mu_1}{\beta}+ \frac{\mu_2}{\beta^2} + 
\cdots +\frac{\mu_n}{\beta^n}, \qquad 
r(\mu) : = \frac{\tilde{\mu}_1}{\beta}+ \frac{\tilde{\mu}_2}{\beta^2} + 
\cdots +\frac{\tilde{\mu}_n}{\beta^n}.
\end{equation*}
If $\tilde{\mu} =\emptyset$,
we set $r(\mu) =1$.
For 
$\mu = (\mu_1,\dots,\mu_n) \in B_n(X_\beta)$,
we set
$s_\mu =s_{\mu_1}\cdots s_{\mu_n}.$
\begin{lem}
For 
$\mu = (\mu_1,\dots,\mu_n) \in B_n(X_\beta)$,
 we have
 \begin{equation}
 s_\mu s_\mu^* = 
 \chi_{[l(\mu),r(\mu))}. \label{eqn:smn}
 \end{equation}
  \end{lem}
 \begin{pf}
For $n=1$, the equality \eqref{eqn:smn} holds by the above proposition.
Suppose that   
the equality \eqref{eqn:smn} holds 
for a fixed $n=k$.
It then follows that 
for $j=0,\dots,N-1$ and $\xi,\eta \in H$, 
\begin{equation}
\langle 
 s_j s_{\mu_1}\cdots s_{\mu_k}s_{\mu_k}^*\cdots s_{\mu_1}^* s_j^* \xi \mid \eta \rangle 
= 
\frac{1}{\beta} 
\int_0^1
\chi_{[l(\mu),r(\mu))}
  \xi(g_j(x)) \overline{\eta(g_j(x))}dx. \label{eqn:int}
\end{equation}
For $j=0,1,\dots,N-2$,
put
$y = g_j(x) \in [\frac{j}{\beta}, \frac{j+1}{\beta}]$ so that
$x = f_\beta(y) = \beta y -j$.
The above equation \eqref{eqn:int} goes to
\begin{equation*}
   \int_0^1 
\chi_{[\frac{j}{\beta},\frac{j+1}{\beta})}(y)
\chi_{[l(\mu),r(\mu))}(f_\beta(y))
 \xi(y) \overline{\eta(y)}dy 
=
\langle 
 \chi_{[\frac{j}{\beta},\frac{j+1}{\beta})\cap f_\beta^{-1}(
[l(\mu),r(\mu)))}
 \xi \mid \eta \rangle. 
\end{equation*}
As
\begin{align*}
& [\frac{j}{\beta},\frac{j+1}{\beta}) \cap f_\beta^{-1}([l(\mu),r(\mu))) \\
=
&[\frac{j}{\beta}+ \frac{\mu_1}{\beta^2} + \frac{\mu_2}{\beta^3}+
\dots +\frac{\mu_k}{\beta^{k+1}},
\frac{j}{\beta}+ \frac{\tilde{\mu}_1}{\beta^2} + \frac{\tilde{\mu}_2}{\beta^3}+
\dots +\frac{\tilde{\mu}_k}{\beta^{k+1}}),
\end{align*}
we have
\begin{equation*}
 s_j s_{\mu_1}\cdots s_{\mu_k}s_{\mu_k}^*\cdots s_{\mu_1}^* s_j^*  
=
\chi_{
[\frac{j}{\beta}+ \frac{\mu_1}{\beta^2} + \frac{\mu_2}{\beta^3}+
\dots +\frac{\mu_k}{\beta^{k+1}},
\frac{j}{\beta}+ \frac{\tilde{\mu}_1}{\beta^2} + 
\frac{\tilde{\mu}_2}{\beta^3}+
\dots +\frac{\tilde{\mu}_k}{\beta^{k+1}})}.
\end{equation*}
Since
the word
$(j,\tilde{\mu}_1, \dots, \tilde{\mu}_k)$
is the minimal word in $B_{k+1}(X_\beta)$
satisfying
$(j,\mu_1,\dots,\mu_k) \prec
(j,\tilde{\mu}_1, \dots, \tilde{\mu}_k)$,
the desired equality holds for $k+1$ and $j=0,\dots,N-2$.
For $j=N-1$, 
 we may similarly show the equality \eqref{eqn:smn}. 
\end{pf}
The following lemma is straightforward.
\begin{lem}
For a measurable subset $F \subset [0,1]$, we have
$s_j^* \chi_F s_j = \chi_{g_j^{-1}(F)}
$ for $j=0,1,\dots,N-1$.
\end{lem}
We then have
\begin{lem}
For the maximal element
$\xi_\beta = (\xi_1,\xi_2, \dots) \in X_\beta$,
we have
\begin{equation}
s_{\xi_1 \xi_2 \cdots \xi_n}^* 
s_{\xi_1 \xi_2 \cdots \xi_n}
=
\chi_{[0,\beta_n]},
\qquad n \in {\mathbb{N}}. \label{eqn:sxis}
\end{equation}
\end{lem}
\begin{pf}
The  equality \eqref{eqn:sxis} holds for $n=1$.
Suppose that the  equality \eqref{eqn:sxis} holds for $n=k$.
It then follows that
\begin{equation*}
s_{\xi_1 \xi_2 \cdots \xi_{k+1}}^* 
s_{\xi_1 \xi_2 \cdots \xi_{k+1}}
=
s_{\xi_{k+1}}^*
\chi_{[0,\beta_k]}
 s_{\xi_{k+1}} 
=
\chi_{g_{\xi_{k+1}}^{-1}([0,\beta_k])}.
\end{equation*}
Since
\begin{equation*}
{g_{\xi_{k+1}}^{-1}([0,\beta_k])}
=
\{x \in [0,1] \mid \frac{1}{\beta} x + \frac{\xi_{k+1}}{\beta} 
\le  \beta^k - \xi_1 \beta^{k-1} - \cdots - \xi_k \}
=[0,\beta_{k+1}],
\end{equation*}
the  equality \eqref{eqn:sxis} holds for $n=k+1$.
\end{pf}

\begin{lem}
For $n \in {\mathbb{N}}$ and 
$j=0,1,\cdots, N-1$,
we have
\begin{equation}
s_{\xi_1 \xi_2 \cdots \xi_n j}^* 
s_{\xi_1 \xi_2 \cdots \xi_n j}
=
\begin{cases}
0 & \text{ for }  j > \xi_{n+1},\\
\chi_{[0,\beta_{n+1}]} & \text{ for } j= \xi_{n+1},\\
1 & \text{ for }  j < \xi_{n+1}.
\end{cases}
\end{equation}
\end{lem}
\begin{pf}
We have
\begin{equation*}
s_{\xi_1 \xi_2 \cdots \xi_n j}^* 
s_{\xi_1 \xi_2 \cdots \xi_n j}
=
s_{\xi_j}^*
\chi_{[0,\beta_n]}
 s_{\xi_j} 
=
\chi_{g_{j}^{-1}([0,\beta_n])}
\end{equation*}
and
\begin{equation*}
{g_{j}^{-1}([0,\beta_n])}
=
\{x \in [0,1] \mid \frac{1}{\beta} x + \frac{j}{\beta} 
\le  \beta_n \}
=[0,\beta \beta_n -j].
\end{equation*}
Since
$\beta \beta_n -j = \beta_{n+1}+\xi_{n+1}-j$,
we have
\begin{equation*}
s_{\xi_1 \xi_2 \cdots \xi_n j}^* 
s_{\xi_1 \xi_2 \cdots \xi_n j}
=
\chi_{[0,\beta_{n+1}+\xi_{n+1}-j]}.
\end{equation*}
If  $\xi_{n+1} = j$,
the equality  
$s_{\xi_1 \xi_2 \cdots \xi_n j}^* 
s_{\xi_1 \xi_2 \cdots \xi_n j}
=
\chi_{[0,\beta_{n+1}]}
$
holds.
If  $\xi_{n+1} < j$,
we have
$\xi_{n+1}- j \le -1$ and hence
$\beta_{n+1} + \xi_{n+1}-j \le 0$
so that
$[0,\beta_{n+1}+\xi_{n+1}-j] = \{0 \}$ or $\emptyset$,
which shows  
$s_{\xi_1 \xi_2 \cdots \xi_n j}^* 
s_{\xi_1 \xi_2 \cdots \xi_n j}
= 0$.
If  $\xi_{n+1} > j$,
we have
$\xi_{n+1}- j \ge 1$ and hence
$\beta_{n+1} + \xi_{n+1}-j \ge 1$
so that
$[0,\beta_{n+1}+\xi_{n+1}-j] = [0,1]$,
which shows  
$s_{\xi_1 \xi_2 \cdots \xi_n j}^* 
s_{\xi_1 \xi_2 \cdots \xi_n j}
= \chi_{[0,1]} =1$.
\end{pf}
Therefore we have
\begin{thm}
The correspondence 
$S_j \longrightarrow s_j$ for $j=0,1,\dots,N-1$
gives rise to an isomorphism from
$\OB$ to the $C^*$-algebra 
$C^*(s_0,s_1,\dots,s_{N-1})$
on $L^2([0,1])$
generated by 
the partial isometries
$s_0,s_1,\dots,s_{N-1}$.
\end{thm}
\begin{pf}
Let us denote by 
${\mathcal{A}}_{[0,1],l}$
the 
$C^*$-algebra on $L^2([0,1])$
generated by the projections
$s_\mu^* s_\mu,
\mu \in B_l(X_\beta)$,
and
${\mathcal{A}}_{[0,1],\beta}$
the 
$C^*$-algebra generated by 
$\cup_{l \in {\mathbb{N}}} {\mathcal{A}}_{[0,1],l}$.
By the previous lemma
and \cite[Corollary 3.2]{KMW},
the $C^*$-algebra
${\mathcal{A}}_{[0,1],l}$
is isomorphic to the $C^*$-subalgebra
${\mathcal{A}}_{l}$ 
of $\OB$
so that
${\mathcal{A}}_{[0,1],\beta}$
is isomorphic to
$\AB$
through the correspondence:
$S_\mu^* S_\mu \longleftrightarrow s_\mu^* s_\mu$ 
for
$\mu \in B_*(X_\beta)$.
The isomorphism from $\AB$ to 
${\mathcal{A}}_{[0,1],\beta}$ is denoted by
$\pi$.
Put
$\rho_j(x) = S_j^* x S_j$ for $x \in \AB, j=0,1,\dots,N-1$.
Then the relations
\begin{equation}
\pi(\rho_j(x)) = s_j^* \pi(x) s_j, \qquad x \in \AB, j=0,1,\dots,N-1
\label{eqn:pirho}
\end{equation}
hold by the previous lemma.
Since the $C^*$-algebra $\OB$ has the universal property 
subject to the relation \eqref{eqn:pirho}
(cf. \cite{MaDocMath02}),
there exists a surjective $*$-homomorphism
$\tilde{\pi}$ from $\OB$
to
$C^*(s_0,s_1,\dots,s_{N-1})$
such that
$\tilde{\pi}(S_j) = s_j, j=0,1,\dots,N-1$
and
$\tilde{\pi}(x) = \pi(x), x \in \AB$.
As the $C^*$-algebra $\OB$ is simple,
the $*$-homomorphism 
$\tilde{\pi}$ is actually an isomorphism.
\end{pf}
In what follows,
we may identify the $C^*$-algebra $\OB$ 
with the $C^*$-algebra 
$C^*(s_0,s_1,\dots,s_{N-1})$
through the identification of the generating partial isometries 
$S_j$ and $s_j, j=0,1,\dots,N-1$.

\section{PL functions for SFT $\beta$-shifts}

In this section,
we will realize the group $\Gamma_\beta$ for an SFT $\beta$-shift
as PL functions on $[0,1)$.
For a word $\mu = (\mu_1, \dots, \mu_n) \in B_n(X_\beta)$,
denote by 
$U_\mu \subset X_\beta$  
the cylinder set 
\begin{equation*}
U_\mu = \{ x=(x_i)_{i \in {\mathbb{N}}} 
\mid x_1 = \mu_1, \dots, x_n = \mu_n \}.
\end{equation*}
We put
\begin{equation*}
\Gamma^+(\mu) = 
\{  (x_i)_{i\in {\mathbb{N}}} \in X_\beta
\mid 
(\mu_1,\dots,\mu_n, x_1, x_2, \dots ) \in X_\beta \}
\end{equation*}
the set of followers of $\mu$.
Recall that $\varphi$ stands for the unique KMS state 
for the gauge action on 
the $C^*$-algebra $\OB$.
We note that the value 
$\varphi(a_{\mu_1\cdots\mu_n})$
is computed inductively in the following way.
For $n=1$, we see
\begin{equation*}
\varphi(a_{\mu_1}) =
\begin{cases}
1 & \text{ if } \mu_1 < \xi_1, \\
\beta - \xi_1 & \text{ if } \mu_1 = \xi_1, \\
0 & \text{ if } \mu_1 > \xi_1.
\end{cases}
\end{equation*}
Suppose that
the value
$\varphi(a_{\mu_1\cdots\mu_k})$
is computed for all
$\mu = (\mu_1,\dots, \mu_k) \in B_k(X_\beta)$
with $ k <n$.
If
 $(\mu_1, \dots, \mu_n)$
is the  maximal element 
$(\xi_1,\dots, \xi_n)$
in $B_n(X_\beta)$,
then 
\begin{equation}
\varphi(a_{\mu_1\cdots\mu_n})
= \beta^n - \xi_1 \beta^{n-1} -\cdots - \xi_{n-1} \beta - \xi_n. \label{eqn:varphiamu}
\end{equation}
If
$(\mu_1, \dots, \mu_n) \ne (\xi_1,\dots, \xi_n)$,
then there exists $k \le n$ such that
$\mu_k < \xi_k$.
If $k=n$, then
 $\varphi(a_{\mu_1\cdots\mu_n})
=1$.
If $k < n$, we see 
$a_{\mu_1\cdots\mu_k} =1$
so that
\begin{equation*}
a_{\mu_1\cdots\mu_n}
= S_{\mu_n}^* \cdots S_{\mu_{k+1}}^* S_{\mu_{k+1}} \cdots S_{\mu_n} =
a_{\mu_{k+1}\cdots\mu_n}. 
\end{equation*}
Hence
\begin{equation*}
\varphi(a_{\mu_1\cdots\mu_n})
= \varphi(
a_{\mu_{k+1}\cdots\mu_n}). 
\end{equation*}
Since $|(\mu_{k+1}, \dots, \mu_n) | <n$,
the value  
$\varphi(a_{\mu_{k+1}\cdots\mu_n})$
is computed.
Therefore the value 
$\varphi(a_{\mu_{1}\cdots\mu_n})$
is computed for all $(\mu_1,\dots, \mu_n) \in B_n(X_\beta)$.
The following lemma is clear
from Lemma 4.5 and \eqref{eqn:varphiamu}.
\begin{lem}
Assume that the generating partial isometries 
$S_0, S_1, \dots, S_{N-1}$
are represented on $L^2([0,1])$.
For a word $\mu \in B_*(X_\beta)$,
 the projection $S_\mu^*S_\mu$ is identified with 
the characteristic function 
$\chi_{[0,\varphi(a_\mu))}$
of the interval $[0,\varphi(a_\mu))$.
\end{lem}
Recall that for a word 
$\nu = (\nu_1,\dots,\nu_n) \in B_n(X_\beta)$,
the  notations 
\begin{equation*}
l(\nu) = \frac{\nu_1}{\beta}+ \frac{\nu_2}{\beta^2} + 
\cdots +\frac{\nu_n}{\beta^n}, \qquad 
r(\nu) = \frac{\tilde{\nu}_1}{\beta}+ \frac{\tilde{\nu}_2}{\beta^2} + 
\cdots +\frac{\tilde{\nu}_n}{\beta^n}
\end{equation*}
are introduced in Section 4
where 
$\tilde{\nu}= (\tilde{\nu}_1,\dots, \tilde{\nu}_n)$
is the smallest word in $B_n(X_\beta)$
satisfying
$\nu \prec\tilde{\nu}$.
If $\nu$ is the maximum word in $B_n(X_\beta)$,
we set $r(\nu) =1$.
The following two lemmas are crucial.
\begin{lem}
For 
$\mu, \nu \in B_*(X_\beta)$,
we have
$\Gamma^+(\mu) = \Gamma^+(\nu)$
if and only if
\begin{equation}
\frac{r(\mu) -l(\mu)}{r(\nu) -l(\nu)} = 
\beta^{|\nu| -|\mu|}.
\end{equation}
\end{lem}
\begin{pf}
We note that 
$\Gamma^+(\mu) = \Gamma^+(\nu)$
if and only if
$S_\mu^* S_\mu = S_\nu^* S_\nu$.
By the above lemma,
we have
$\Gamma^+(\mu) = \Gamma^+(\nu)$
if and only if
$\varphi(a_\mu) =\varphi(a_\nu)$.
Since
$
\varphi(S_\mu^*S_\mu) = \beta^{|\mu|}\varphi( S_\mu S_\mu^*)
$
and
$\varphi( S_\mu S_\mu^*) = r(\mu) -l(\mu)$,
we have
$\varphi(a_\mu) =\varphi(a_\nu)$
if and only if
$\beta^{|\mu|}
(r(\mu) -l(\mu))
=
\beta^{|\nu|}
(r(\nu) -l(\nu)).
$
\end{pf}
We note that the above lemma holds for any real number $\beta >1$
even if $(X_\beta,\sigma)$ is not a shift of finite type.
\begin{lem}\label{lem:ut}
For $\tau \in \Gamma_\beta$, 
there exists $u_\tau \in N(\DB,\OB)$
such that 
there exist
$\mu(i), \nu(i) \in B_*(X_\beta), i=1,\dots,m$  
satisfying
\begin{enumerate}
\renewcommand{\labelenumi}{(\arabic{enumi})}
\item
$ u_\tau = \sum_{i=1}^m S_{\mu(i)}S_{\nu(i)}^*$
such that
{\begin{enumerate}
\renewcommand{\labelenumi}{(\alpha{enumi})}
\item
$S_{\nu(i)}^*S_{\nu(i)}
 = S_{\mu(i)}^*S_{\mu(i)}, \quad i=1,\dots,m,
$
 \item
$\sum_{i=1}^m S_{\nu(i)}S_{\nu(i)}^* 
 =\sum_{i=1}^m S_{\mu(i)}S_{\mu(i)}^*=1.
$
\end{enumerate}}
\item
$f \circ \tau^{-1} = u_\tau f u_\tau^*$ for $f \in \DB$.
\end{enumerate}
\end{lem}
\begin{pf}
Since
$(X_\beta,\sigma)$ is SFT,
there exist continuous functions
$k, l: X_\beta \longrightarrow \Zp$
for $\tau \in \Gamma_\beta$
such that
$\sigma^{l(x)}(\tau(x)) = \sigma^{k(x)}(x), x \in X_\beta$. 
Hence there exists a family of cylinder sets
$U_{\nu(1)}, \dots, U_{\nu(m)}, U_{\mu(1)}, \dots, U_{\mu(m)}$
such that
\begin{gather*}
\Gamma^+(\nu(i))  = \Gamma^+(\mu(i)), \qquad i=1,\dots,m, \\
X_\beta  = \sqcup_{i=1}^m U_{\nu(i)}
         = \sqcup_{i=1}^m U_{\mu(i)}
\end{gather*}
and
\begin{equation*}
\tau(x_1,x_2,\dots )
=(\mu(i)_1,\dots,\mu(i)_{l_i}, x_{k_i+1}, x_{k_i+2}, \dots )
\text{ for }
(x_n)_{n \in {\mathbb{N}}} \in U_{\nu(i)}
\end{equation*}
where
$l_i = |\mu(i)|, k_i =|\nu(i)|$
and
$\mu(i) = (\mu(i)_1,\dots,\mu(i)_{l_i})$.
Hence we have
\begin{equation*}
\sum_{i=1}^m S_{\nu(i)}S_{\nu(i)}^* =
\sum_{i=1}^m S_{\mu(i)}S_{\mu(i)}^*=1,
\qquad 
S_{\nu(i)}^*S_{\nu(i)} = S_{\mu(i)}^*S_{\mu(i)}, \quad i=1,\dots,m.
\end{equation*}
By putting
$ u_\tau = \sum_{i=1}^m S_{\mu(i)}S_{\nu(i)}^*$
we see that $u_\tau$ belongs to
$N(\DB,\OB)$
and
satisfies
$\chi_{U_\eta} \circ\tau^{-1}
= u_\tau S_\eta S_\eta^* u_\tau^*$
for all $\eta \in B_*(X_\beta)$
so that
$f \circ \tau^{-1} = u_\tau f u_\tau^*$ for $f \in \DB$.
\end{pf}

Following Nekrashevych in \cite{Nek},
we will introduce 
a notation of tables in order to  represent elements of $\Gamma_\beta$.

\noindent
{\bf Definition.}
A $\beta$-{\it adic table} for
SFT $\beta$-shift is a matrix
\begin{equation*}
\begin{bmatrix}
\mu(1) & \mu(2) & \cdots & \mu(m) \\
\nu(1) & \nu(2) & \cdots & \nu(m) 
\end{bmatrix}
\end{equation*}
for 
$\nu(i), \mu(i) \in B_*(X_\beta), i=1,\dots, m$
such that
\begin{enumerate}
\renewcommand{\labelenumi}{(\alph{enumi})}
\item $\Gamma^+(\nu(i)) = \Gamma^+(\mu(i)), i=1,\dots,m$,
\item $X_\beta = \sqcup_{i=1}^m U_{\nu(i)}
               = \sqcup_{i=1}^m U_{\mu(i)} :$
               disjoint unions.
\end{enumerate}
We may assume that 
$\nu(1) \prec \nu(2) \prec \cdots \prec  \nu(m)$.
Since the above two conditions (a), (b) are equivalent to the conditions (a), (b) in
Lemma \ref{lem:ut} (1) respectively, 
we have
\begin{lem}
For an element $\tau \in \Gamma_\beta$ with its unitary
$ u_\tau = \sum_{i=1}^m S_{\mu(i)}S_{\nu(i)}^* \in N(\DB,\OB)$
as in Lemma \ref{lem:ut}, 
the matrix
\begin{equation*}
T_\tau =
\begin{bmatrix}
\mu(1) & \mu(2) & \cdots & \mu(m) \\
\nu(1) & \nu(2) & \cdots & \nu(m) 
\end{bmatrix}
\end{equation*}
is a $\beta$-adic table for SFT $\beta$-shift. 
\end{lem}

\noindent
{\bf Definition.}
\begin{enumerate}
\renewcommand{\labelenumi}{(\roman{enumi})}
\item 
An interval $[x_1,x_2)$ in $[0,1]$ is said to be 
a $\beta$-{\it adic interval} for word $\nu \in B_*(X_\beta)$ 
if 
$x_1 =l(\nu) $ and
$x_2 = r(\nu)$.
\item
A rectangle
$I \times J$ in 
$[0,1]\times [0,1]$
 is said to be 
a $\beta$-{\it adic rectangle}
if  both 
$I, J$ 
are $\beta$-adic intervals for 
words
$\nu \in B_n(X_\beta), \mu \in B_m(X_\beta)$
such that
$I = [l(\nu), r(\nu))$,
$J = [l(\mu), r(\mu))$
and
\begin{equation*}
\frac{r(\mu) - l(\mu)}{r(\nu) - l(\nu)} = \beta^{n -m} 
\end{equation*}
\item
For two partitions
$0=x_0 < x_1 < \dots <x_{m-1} < x_m = 1$
and
$0=y_0 < y_1 < \dots <y_{m-1} < y_m = 1$
of $[0,1]$,
put
$I_p = [x_{p-1},x_p),
J_p = [y_{p-1},y_p),
p=1,2,\dots,m$.
The partition
$I_p \times J_q, p,q =1,\dots,m$ of 
$[0,1)\times [0,1)$
is said to be 
a $\beta$-{\it adic pattern of squares } 
for SFT $\beta$-shift
if there exists a permutation 
$\sigma $ on $\{1,2,\dots,m\}$
such that
the rectangles
$I_p \times J_{\sigma(p)}$
are $\beta$-adic rectangles for all $p=1,2,\dots,m$.
\end{enumerate}
For a 
$\beta$-adic pattern of squares above,
the slopes of diagonals 
$s_p = \frac{y_{\sigma(p)} -y_{\sigma(p)-1}}{x_p -x_{p-1}}, p=1,2,\dots,m$
are said to be rectangle slopes.
We then have
\begin{lem}\label{lem:tableSFT}
For a $\beta$-adic table 
\begin{equation*}
T =
\begin{bmatrix}
\mu(1) & \mu(2) & \cdots & \mu(m) \\
\nu(1) & \nu(2) & \cdots & \nu(m) 
\end{bmatrix},
\end{equation*}
there exists a $\beta$-adic pattern of squares 
whose rectangle slopes are
\begin{equation*}
\beta^{|\nu(1)| -|\mu(1)|}, 
\beta^{|\nu(2)| -|\mu(2)|},
\dots, 
\beta^{|\nu(m)| -|\mu(m)|}
\end{equation*}
\end{lem}
\begin{pf}
We are assuming the ordering such as
$\nu(1) \prec\dots \prec\nu(m)$.
Since
$X_\beta = \sqcup_{j=1}^m U_{\mu(j)}$ disjoint union,
there exists a permutation $\sigma_0$ on $\{1,2,\dots,m\}$
such that
$
\mu(\sigma_0(1)) \prec 
\mu(\sigma_0(2)) \prec
\dots \prec
\mu(\sigma_0(m))$.
Put
\begin{equation*}
x_i = l(\nu(i+1)), 
\qquad 
y_i = l(\mu(\sigma_0(i+1))), \quad i= 0, 1,\dots,m-1
\end{equation*}
and
\begin{equation*}
I_p =[x_{p-1}, x_p), \qquad 
J_p =[y_{p-1}, y_p), 
 \quad p=1,2,\dots,m.
\end{equation*}
Define the permutation 
$\sigma := \sigma_0^{-1}$ on $\{1,2,\dots,m\}$.
We note that
$r(\nu(i)) = l(\nu(i+1)),
r(\mu(\sigma_0(i))) = l(\mu(\sigma_0(i+1)))
$
for
$
i=1,\dots,m-1$.
Then the rectangles
$I_p \times J_{\sigma(p)}, p=1,2,\dots,m$
are $\beta$-adic rectangles such that
\begin{equation*}
\frac{y_{\sigma(p)} -y_{\sigma(p)-1}}{x_p -x_{p-1}}
=
\frac{r(\mu(p)) -l(\mu(p))}{r(\nu(p)) -l(\nu(p))}.
\end{equation*}
Since
$
r(\zeta) - l(\zeta) 
= \varphi(S_\zeta S_\zeta^*)
= \frac{1}{\beta^{ |\zeta|}} \varphi(S_\zeta^*S_\zeta)
$
for $\zeta \in B_*(X_\beta)$,
we have
\begin{align*}
r(\nu(p)) - l(\nu(p)) 
& = \frac{1}{\beta^{ |\nu(p)|}} \varphi(S_{\nu(p)}^*S_{\nu(p)}),\\
r(\mu(p)) - l(\mu(p)) 
& = \frac{1}{\beta^{|\mu(p)|}} \varphi(S_{\mu(p)}^*S_{\mu(p)}).
\end{align*}
As the condition
$\Gamma^+(\nu(p)) = \Gamma^+(\mu(p))$ 
implies 
$S_{\nu(p)}^*S_{\nu(p)}= S_{\mu(p)}^*S_{\mu(p)}$,
we have
\begin{equation*}
\frac{y_{\sigma(p)} -y_{\sigma(p)-1}}{x_p -x_{p-1}}
=
\beta^{|\nu(p)| - |\mu(p)|},
\qquad p=1,2,\dots,m.
\end{equation*}
\end{pf}
We define a $\beta$-adic version of piecewise linear functions on $[0,1)$
in the following way.

\noindent
{\bf Definition.}
A piecewise linear function $f$ on $[0,1)$
is called a $\beta$-adic PL function for SFT $\beta$-shift
if $f$ is a right continuous bijection on $[0,1)$
such that there exists a $\beta$-adic pattern of squares 
$I_p \times J_p, p=1,2,\dots,m$
where
$I_p = [x_{p-1}, x_p), J_p = [y_{p-1},y_p), p=1,\dots,m$ 
with a permutation $\sigma$ on $\{1,2,\dots,m \}$
such that
\begin{equation*}
f(x_{p-1}) = y_{\sigma(p)-1},
\qquad
f_{-}(x_p) = y_{\sigma(p-1)+1},
\qquad
 p=1,2,\dots,m
\end{equation*}
where
$
f_{-}(x_p) = \lim_{h \to 0+} f(x_p -h),
$
and $f$ is linear on $[x_{p-1}, x_p)$
with slope
$\frac{y_{\sigma(p)} -y_{\sigma(p)-1}}{x_p -x_{p-1}}$
for
$ p=1,2,\dots,m$.

The following proposition is immediate from 
the definition of $\beta$-adic PL functions.
\begin{prop}
A $\beta$-adic PL function for SFT $\beta$-shift naturally 
gives rise to
a $\beta$-adic pattern of squares for SFT $\beta$-shift.
\end{prop}
We may directly construct a $\beta$-adic PL function $f_T$ from a $\beta$-adic table 
\begin{math}
T=
\bigl[
\begin{smallmatrix}
\mu(1) & \mu(2) & \cdots & \mu(m) \\
\nu(1) & \nu(2) & \cdots & \nu(m) 
\end{smallmatrix}
\bigr]
\end{math}
as follows.
Put
$x_i = l(\nu(i+1)), \hat{y}_i = l(\mu(i+1)), i= 0,1,\dots, m-1$.
Define $f_T$ by
$f_T(x_i) = \hat{y}_i, i=0, 1,\dots,m-1$
and
$f_T$ is linear on $[x_{i-1},x_i), i=1,2,\dots,m$
with slope
$\frac{\hat{y}_i - \hat{y}_{i-1}}{x_i - x_{i-1}} 
=\frac{r(\mu(i)) - l(\mu(i))}{r(\nu(i)) - l(\nu(i))} 
= \beta^{|\nu(i)| -|\mu(i)|}.
$
Hence the function $f_T$ yields a $\beta$-adic PL function.

It is straightforward to see that 
the composition of two $\beta$-adic PL functions
is also a $\beta$-adic PL function.
Hence
the set of $\beta$-adic PL functions forms a group
under compositions.
We reach the following:
\begin{thm}\label{thm:SFTPL}
The topological full group
$\Gamma_\beta$
for a SFT $\beta$-shift $(X_\beta,\sigma)$
is realized as the group of all $\beta$-adic PL functions
for SFT $\beta$-shift.
\end{thm}

\section{PL functions for sofic $\beta$-shifts}
In this section, we will represent the topological full group $\Gamma_\beta$
for sofic $\beta$-shifts as piecewise linear functions on $[0,1)$.
Througout this section, we assume that $(X_\beta,\sigma)$ is sofic.
By Lemma 2.2, the algebra $\AB$ is finite dimensional.
We set $K_\beta = \dim \AB$.
Let $E_1,\dots,E_{K_\beta}$
be the minimal projections of $\AB$ so that
$\sum_{i=1}^{K_\beta} E_i = 1$.
Then any minimal projection
$E_i$ is of the form
$E_i = a_{\xi_1\cdots\xi_{p_i}} - a_{\xi_1\cdots\xi_{q_i}}$
for some
$p_i, q_i \in \Zp$.
We order 
$E_1, \dots,E_{K_\beta}$ 
following the order
$
   \varphi(a_{\xi_1\cdots\xi_{p_1}}) 
< \cdots
< \varphi(a_{\xi_1\cdots\xi_{p_{K_\beta}}})$
in ${\mathbb{R}}$,
where $\varphi$ is the unique KMS-state on $\OB$ for the gauge action.
Recall that $\hat{\rho}_t \in \Aut(\OB), t \in {\mathbb{R}}/{\mathbb{Z}}$
denotes the gauge action on $\OB$,
and
$N(\DB,\OB)$ denotes the normalizer group of $\DB \subset \OB$.
Fix $v\in N(\DB,\OB)$ for a while.
For $m \in {\mathbb{Z}}$ 
and $\mu \in B_n(X_\beta), n \in {\mathbb{N}}$, 
put
\begin{equation*}
v_m  = \int_{\Bbb T}\hat{\rho}_t(v) e^{-2\pi\sqrt{-1}mt} dt 
\quad
\text{ and }  
\quad
v_\mu  = S_\mu^* v_n, 
\quad v_{-\mu} =  v_{-n}S_\mu. 
\end{equation*}
It is straightforward to see the following lemma.
\begin{lem}
The operators 
$ v_\mu , v_{-\mu}$
for $\mu \in B_n(X_\beta)$
and
$v_0$ 
are partial isometries in $\FB$
such that $v$ is decomposed as a following finite sum:
\begin{equation*}
v = 
\sum_{n finite}\sum_{\mu \in B_n(X_\beta)} S_\mu v_\mu
+ 
v_0
+
\sum_{n finite}\sum_{\mu \in B_n(X_\beta)} v_{-\mu}S_\mu^*
\end{equation*}
such that
$v_\mu \DB v_\mu^*, v_\mu^* \DB v_\mu, v_{-\mu} \DB v_{-\mu}^*$ 
and $v_{-\mu}^* \DB v_{-\mu}$ are contained in $\DB$. 
\end{lem}
Define a subalgebra ${\mathcal{F}}_\beta^k$ of $\FB$ for $k \in \Zp$
by 
\begin{equation*}
{\mathcal{F}}_\beta^k =
C^*(S_\xi E_i S_\eta^* \mid \xi, \eta \in B_k(X_\beta), i=1,\dots,K_\beta). 
\end{equation*}
We set
\begin{equation*}
\supp_+(v)  = \{ \mu \in B_*(X_\beta) \mid v_\mu \ne 0 \},\qquad
\supp_-(v)  = \{ \mu \in B_*(X_\beta) \mid v_{-\mu} \ne 0 \}.
\end{equation*}
Both of them are finite sets.
For $\mu \in \supp_+(v)$,
there exists $k_+(\mu) \in \Zp$ such that
$v_\mu \in {\mathcal{F}}_\beta^{k_+(\mu)}$.
For $\mu \in \supp_-(v)$,
there exists $k_-(\mu) \in \Zp$ such that
$v_{-\mu} \in {\mathcal{F}}_\beta^{k_-(\mu)}$.
There exists $k_0 \in \Zp$ such that
$v_0 \in {\mathcal{F}}_\beta^{k_0}$.
We then have
\begin{lem} Keep the above notations.
\begin{enumerate}
\renewcommand{\labelenumi}{(\roman{enumi})}
\item
For $\mu \in \supp_+(v)$
and
$\eta \in B_{k_+(\mu)}(X_\beta), i=1,\dots,K_\beta$
such that
$v_\mu^*v_\mu 
\ge 
S_\eta E_i S_\eta^* \ne 0,
$ 
there uniquely exists
$\xi \in B_{k_+(\mu)}(X_\beta)$
such that
$v_\mu v_\mu^* 
\ge 
S_\xi E_i S_\xi^* \ne 0
$ 
and
\begin{equation*}
Ad(v_\mu) 
(S_\eta E_i S_\eta^*) =S_\xi E_i S_\xi^*.
\end{equation*}
\item
For $\mu \in \supp_-(v)$
and
$\eta \in B_{k_-(\mu)}(X_\beta), i=1,\dots,K_\beta$
such that
$v_{-\mu}^*v_{-\mu} 
\ge 
S_\eta E_i S_\eta^*\ne 0,
$ 
there uniquely exists
$\xi \in B_{k_-(\mu)}(X_\beta)$
such that
$v_{-\mu} v_{-\mu}^* 
\ge 
S_\xi E_i S_\xi^*\ne 0
$ 
and
\begin{equation*}
Ad(v_{-\mu}) 
(S_\eta E_i S_\eta^*) =S_\xi E_i S_\xi^*.
\end{equation*}
\item
For
$\eta \in B_{k_0}(X_\beta), i=1,\dots,K_\beta$
such that
$v_{0}^*v_{0} 
\ge 
S_\eta E_i S_\eta^*\ne 0,
$ 
there uniquely exists
$\xi \in B_{k_0}(X_\beta)$
such that
$v_{0} v_{0}^* 
\ge 
S_\xi E_i S_\xi^*\ne 0
$ 
and
\begin{equation*}
Ad(v_{0}) 
(S_\eta E_i S_\eta^*) =S_\xi E_i S_\xi^*.
\end{equation*}
\end{enumerate}
\end{lem}
\begin{pf}
(i)
As $v_\mu \in {\mathcal{F}}_\beta^{k_+(\mu)}$,
it is written
$
v_\mu = \sum_{\xi,\eta' \in B_{k_+(\mu)}(X_\beta)}S_\xi a_{\xi,\eta'}S_{\eta'}^*
$
for some
$a_{\xi,\eta'} \in \AB$.
Suppose that
$v_\mu^*v_\mu 
\ge 
S_\eta E_i S_\eta^* \ne 0.
$ 
Hence 
$
S_\eta^* S_\eta \ge E_i.
$ 
It then follows that
\begin{align*}
Ad(v_\mu) 
(S_\eta E_i S_\eta^*) 
& = v_\mu S_\eta E_i S_\eta^* v_\mu^* \\
& = \sum_{\xi,\xi' \in B_{k_+(\mu)}(X_\beta)}S_\xi a_{\xi,\eta}S_{\eta}^*
S_\eta E_i S_\eta^* S_\eta a_{\xi',\eta}^* S_{\xi'}^* \\
& = \sum_{\xi,\xi' \in B_{k_+(\mu)}(X_\beta)}S_\xi a_{\xi,\eta} E_i 
 a_{\xi',\eta}^* S_{\xi'}^*. 
\end{align*} 
Since 
$Ad(v_\mu)(S_\eta E_i S_\eta^*)$
belongs to  $\DB$,
we have
for $\xi \ne \xi'$ 
\begin{equation*}
0 =
S_\xi S_\xi^*
Ad(v_\mu)(S_\eta E_i S_\eta^*) 
S_{\xi'}S_{\xi'}^*  
=
S_\xi a_{\xi,\eta} E_i a_{\xi',\eta}^* 
 S_{\xi'}^*
\end{equation*} 
so that
\begin{equation*}
Ad(v_\mu) 
(S_\eta E_i S_\eta^*) 
 = \sum_{\xi \in B_{k_+(\mu)}(X_\beta)}S_\xi a_{\xi,\eta} E_i 
 a_{\xi,\eta}^* S_{\xi}^*. 
\end{equation*} 
Since
$v_\mu v_\mu^* 
 = \sum_{\xi,\zeta \in B_{k_+(\mu)}(X_\beta)}
S_\xi a_{\xi,\zeta} a_{\xi,\zeta}^* S_\xi^*
$
is a projection,
the operators
$a_{\xi,\eta} a_{\xi,\eta}^*$ are projections
in $\AB$ for all
$\xi \in B_{k_+(\mu)}(X_\beta)$.
As 
$
S_\xi^*S_\xi a_{\xi,\eta} E_i a_{\xi',\eta}^* 
 S_{\xi'}^*S_{\xi'}
=a_{\xi,\eta} E_i a_{\xi',\eta}^*,
$
we have
$
a_{\xi,\eta} a_{\xi,\eta}^*
\cdot
a_{\xi',\eta} a_{\xi',\eta}^* 
=0
$
for 
$\xi \ne \xi'$
 so that
there uniquely exists
$\xi \in B_{k_+(\mu)}(X_\beta)$
such that
$
a_{\xi,\eta} a_{\xi,\eta}^* E_i =E_i
$
for the word
$\eta$
and $i$.
By the identity
$
a_{\xi,\eta} E_i a_{\xi,\eta}^*
=
a_{\xi,\eta} a_{\xi,\eta}^* E_i,
$
we have
\begin{equation*}
Ad(v_\mu) 
(S_\eta E_i S_\eta^*) 
 = S_\xi E_i S_{\xi}^*. 
\end{equation*} 
(ii),(iii) are similar to (i).
\end{pf}
\begin{prop}
For a  unitary
$v \in N(\DB,\OB)$,
there exists a finite family of  partial isometries
$v_\mu, v_0, v_{-\mu}$ in 
$\FB$ such that $v$ 
is decomposed in the following way:
\begin{equation*}
v = 
\sum_{n finite}\sum_{\mu \in B_n(X_\beta)} S_\mu v_\mu
+ 
v_0
+
\sum_{n finite}\sum_{\mu \in B_n(X_\beta)} v_{-\mu}S_\mu^*
\end{equation*}
such that
\begin{enumerate}
\renewcommand{\labelenumi}{(\arabic{enumi})}
\item 
For any
$\eta \in B_{k_+(\mu)}(X_\beta)$
with $S_\eta E_i S_\eta^* \le v_\mu^* v_\mu$,
the equality
$$
Ad(S_\mu v_\mu)(S_\eta E_i S_\eta^*)= S_\mu S_\xi E_i S_\xi^*S_\mu^*
$$
holds for some 
$\xi \in B_{k_+(\mu)}(X_\beta)$.
\item 
For any
$\eta \in B_{k_0}(X_\beta)$
with 
$S_\eta E_i S_\eta^* \le v_0^* v_0$,
the equality
$$
Ad(v_0)(S_\eta E_i S_\eta^*)= S_\xi E_i S_\xi^*
$$
holds for some 
$\xi \in B_{k_0}(X_\beta)$.
\item 
For any
$\eta \in B_{k_-(\mu)}(X_\beta)$
with 
$S_\eta E_i S_\eta^* \le v_{-\mu}^* v_{-\mu}$,
the equality
$$
Ad(v_{-\mu}S_\mu^*)(S_\mu S_\eta E_i S_\eta^* S_\mu^*)
= S_\xi E_i S_\xi^*
$$
holds
for some 
$\xi \in B_{k_-(\mu)}(X_\beta)$.
\end{enumerate}
\end{prop}
Therefore we have
\begin{lem}\label{lem:utsofic}
For $\tau \in \Gamma_\beta$, 
there exists $v_\tau \in N(\DB,\OB)$
such that 
there exists a family 
$S_{\nu(j)} E_{i_j} S_{\nu(j)}^*,S_{\mu(j)} E_{i_j} S_{\mu(j)}^*, 
j=1,\dots,m$  
of projections satisfying
\begin{enumerate}
\renewcommand{\labelenumi}{(\arabic{enumi})}
\item
$ v_\tau = \sum_{j=1}^m S_{\mu(j)} E_{i_j} S_{\nu(j)}^*$
such that
{\begin{enumerate}
\renewcommand{\labelenumi}{(\alph{enumi})}
\item
$S_{\nu(j)}^*S_{\nu(j)}, S_{\mu(j)}^*S_{\mu(j)} \ge E_{i_j}, \quad j=1,\dots,m,$
\item
$
\sum_{j=1}^m S_{\nu(j)} E_{i_j} S_{\nu(j)}^* =
\sum_{j=1}^m S_{\mu(j)} E_{i_j} S_{\mu(j)}^*=1.
$
\end{enumerate}}
\item
$f \circ \tau^{-1} = v_\tau f v_\tau^*$ for $f \in \DB$.
\end{enumerate}
\end{lem}

For $i=1,2,\dots,K_\beta$, put
\begin{equation*}
\Gamma_n^-(i)  = \{ \mu \in B_n(X_\beta) \mid
S_\mu^* S_\mu \ge E_i \},\qquad                    
\Gamma_*^-(i)  = \cup_{n=0}^\infty \Gamma_n^-(i).
\end{equation*}
For $\nu = (\nu_1,\dots, \nu_n) \in \Gamma_n^-(i)$ and
$i=1,\dots,K_\beta$,
put the projection in $\DB$ 
\begin{equation*}
\nu_{[i]}:= S_\nu E_i S_\nu^*
\end{equation*}
and define 
\begin{align*}
r(\nu_{[i]}) 
& = l(\nu) + \frac{1}{\beta^n} \varphi(a_{\xi_1\cdots\xi_{p_i}}) \\
& = \frac{\nu_1}{\beta} + \frac{\nu_2}{\beta^2} +
\cdots + \frac{\nu_n}{\beta^n}
+
 \frac{\xi_{p_i +1}}{\beta^{n+1}} + \frac{\xi_{p_i +2}}{\beta^{n+2}} +
\cdots, \\
l(\nu_{[i]}) 
& = l(\nu) + \frac{1}{\beta^n}\varphi(a_{\xi_1\cdots\xi_{q_i}}) \\
& = \frac{\nu_1}{\beta} + \frac{\nu_2}{\beta^2} +
\cdots + \frac{\nu_n}{\beta^n}
+
 \frac{\xi_{q_i +1}}{\beta^{n+1}} + \frac{\xi_{q_i +2}}{\beta^{n+2}} +
\cdots
\end{align*}
where
$E_i = a_{\xi_1\cdots\xi_{p_i}} -a_{\xi_1\cdots\xi_{q_i}}.$
The following lemma is obvious.
\begin{lem}
Assume that the generating partial isometries 
$S_0, S_1, \cdots, S_{N-1}$
are represented on $L^2([0,1])$.
For $\nu \in \Gamma_n^-(i)$, 
the projection $S_\nu E_i S_\nu^*$ is identified with 
the characteristic function 
$\chi_{[ l(\nu_{[i]}), r(\nu_{[i]}))}$
of the interval $[ l(\nu_{[i]}), r(\nu_{[i]}))$. 
\end{lem}
For $\nu \in \Gamma_*^-(i)$ and
$\mu \in \Gamma_*^-(j)$ with
$S_\nu E_i S_\nu^* \cdot S_\mu E_j S_\mu^*=0$,
define
\begin{equation*}
\nu_{[i]} < \mu_{[j]}  \text{ if } r(\nu_{[i]}) \le l(\mu_{[j]}).
\end{equation*}
Note that
under the condition
$S_\nu E_i S_\nu^* \cdot S_\mu E_j S_\mu^* =0$,
the intervals 
$[l(\nu_{[i]}), r(\nu_{[i]}))$
and
$[l(\mu_{[j]}), r(\mu_{[j]}))$
are disjoint. 
Hence 
the condition 
$\nu_{[i]} < \mu_{[j]}$
implies that
the interval
$[l(\nu_{[i]}), r(\nu_{[i]}))$
is located 
the left side of 
$[l(\mu_{[j]}), r(\mu_{[j]}))$.
\begin{lem}
Keep the above notations.
\begin{enumerate}
\renewcommand{\labelenumi}{(\roman{enumi})}
\item
For $\nu \in \Gamma_n^-(i)$ and
$\mu \in \Gamma_k^-(j)$ we have
$S_\nu E_i S_\nu^* \cdot S_\mu E_jS_\mu^* =0$
if and only if
$[l(\nu_{[i]}), r(\nu_{[i]})) \cap 
[l(\mu_{[j]}), r(\mu_{[j]})) =\emptyset.
$ 
\item
For
$
\nu(j) \in \Gamma_{n_j}^-(i_j), j=1,\dots,m,
$ 
 we have
$ \sum_{j=1}^m S_{\nu(j)} E_{i_j} S_{\nu(j)}^* =1$
if and only if
$
[0,1) =
\sqcup_{j=1}^m
[l(\nu(j)_{[i_j]}), r(\nu(j)_{[i_j]}))
$ a disjoint union.
\item
For
$
\nu(j) \in \Gamma_{n_j}^-(i_j), j=1,\dots,m
$ 
such that
$ \sum_{j=1}^m S_{\nu(j)} E_{i_j} S_{\nu(j)}^* =1$
and
$$
\nu(1)_{[i_1]} < \nu(2)_{[i_2]} < \cdots < \nu(m)_{[i_m]},
$$
we have
$$
r(\nu(j)_{[i_j]}) = l(\nu(j+1)_{[i_{j+1}]}),
\qquad j=1,\dots,m.
$$
\end{enumerate}
\end{lem}

\noindent
{\bf Definition.}
A $\beta$-{\it adic table }for sofic $\beta$-shift is a matrix 
\begin{equation*}
T =
\begin{bmatrix}
{\mu(1)_{[i_1]}} & {\mu(2)_{[i_2]}} & \cdots & {\mu(m)_{[i_m]}} \\
{\nu(1)_{[i_1]}} & {\nu(2)_{[i_2]}} & \cdots & {\nu(m)_{[i_m]}} 
\end{bmatrix}
\end{equation*}
such that
\begin{enumerate}
\renewcommand{\labelenumi}{(\alph{enumi})}
\item
$
\nu(j) \in \Gamma_*^-(i_j), \quad 
\mu(j) \in \Gamma_*^-(i_j) \quad
\text{ for } j=1,\dots,m.
$
\item
$\sqcup_{j=1}^m
[l(\nu(j)_{[i_j]}), r(\nu(j)_{[i_j]})) =
\cup_{j=1}^m
[l(\mu(j)_{[i_j]}), r(\mu(j)_{[i_j]})) =[0,1).
$
\end{enumerate}
We may assume that
\begin{equation*}
\nu(1)_{[i_1]} < \nu(2)_{[i_2]} < \cdots < \nu(m)_{[i_m]}. 
\end{equation*}
Therefore we have 
\begin{lem}
For an element $\tau \in \Gamma_\beta$ with its unitary
$ v_\tau = \sum_{j=1}^m S_{\mu(j)} E_{i_j} S_{\nu(j)}^* \in N(\DB,\OB)$
as in Lemma \ref{lem:utsofic}, 
the matrix
\begin{equation*}
T_\tau =
\begin{bmatrix}
{\mu(1)_{[i_1]}} & {\mu(2)_{[i_2]}} & \cdots & {\mu(m)_{[i_m]}} \\
{\nu(1)_{[i_1]}} & {\nu(2)_{[i_2]}} & \cdots & {\nu(m)_{[i_m]}} 
\end{bmatrix}
\end{equation*}
is a $\beta$-adic table for sofic $\beta$-shift. 
\end{lem}

\noindent
{\bf Definition.}
\begin{enumerate}
\renewcommand{\labelenumi}{(\roman{enumi})}
\item 
An interval $[x_1,x_2)$ in $[0,1]$ is said to be 
a $\beta$-{\it adic interval}
for word $\nu_{[i]}$
if 
$x_1 =l(\nu_{[i]}) $ and
$x_2 = r(\nu_{[i]})$
for some 
$\nu \in B_*(X_\beta)$
and $i =1,\dots,K_\beta$.
\item
A rectangle
$I \times J$ in 
$[0,1]\times [0,1]$
 is said to be 
a $\beta$-{\it adic rectangle}
if  both 
$I, J$ 
are $\beta$-adic intervals for 
words
$\nu_{[i]}, \mu_{[i]}$
such that
$I = [l(\nu_{[i]}), r(\nu_{[i]}))$
and
$J = [l(\mu_{[i]}), r(\mu_{[i]}))$
and
\begin{equation*}
\frac{r(\mu_{[i]}) - l(\mu_{[i]})}{r(\nu_{[i]}) - l(\nu_{[i]})} = \beta^{|\nu|- |\mu|} 
\end{equation*}
\item
For two partitions
$0=x_0 < x_1 < \dots <x_{m-1} < x_m = 1$
and
$0=y_0 < y_1 < \dots <y_{m-1} < y_m = 1$
of $[0,1]$,
put
$I_p = [x_{p-1},x_p),
J_p = [y_{p-1},y_p),
p=1,2,\dots,m$.
The partition
$I_p \times J_q, p,q =1,\dots,m$ of 
$[0,1)\times [0,1)$
is said to be a
$\beta$-{\it adic pattern of squares }
for sofic $\beta$-shift 
if there exists a permutation 
$\sigma $ on $\{1,2,\dots,m\}$
such that
the rectangles
$I_p \times J_{\sigma(p)}$
are $\beta$-adic rectangles for all $p=1,2,\dots,m$.
\end{enumerate}
For a $\beta$-adic pattern of squares above,
the slopes of diagonals 
$s_p = \frac{y_{\sigma(p)}  -y_{\sigma(p)-1}}{x_p -x_{p-1}}, p=1,2,\dots,m$
are said to be rectangle slopes.
Similarly to Lemma \ref{lem:tableSFT} for SFT $\beta$-shift, 
we have
\begin{lem}
For a $\beta$-adic table for sofic $\beta$-shift
\begin{equation*}
T =
\begin{bmatrix}
{\mu(1)_{[i_1]}} & {\mu(2)_{[i_2]}} & \cdots & {\mu(m)_{[i_m]}} \\
{\nu(1)_{[i_1]}} & {\nu(2)_{[i_2]}} & \cdots & {\nu(m)_{[i_m]}} 
\end{bmatrix},
\end{equation*}
there exists a $\beta$-adic pattern of squares for sofic $\beta$-shift
whose rectangle slopes are
\begin{equation*}
\beta^{|\nu(1)| -|\mu(1)|}, 
\beta^{|\nu(2)| -|\mu(2)|},
\dots, 
\beta^{|\nu(m)| -|\mu(m)|}
\end{equation*}
\end{lem}

Similarly to the preceding section,
we will define $\beta$-adic versions of piecewise linear functions on $[0,1)$
for sofic $\beta$-shift
in the following way.

\noindent
{\bf Definition.}
A piecewise linear function $f$ on $[0,1)$
is called a $\beta$-adic PL function for sofic $\beta$-shift
if $f$ is a right continuous bijection on $[0,1)$
such that there exists a $\beta$-adic pattern of squares 
$I_p \times J_p, p=1,2,\dots,m$
where
$I_p = [x_{p-1}, x_p), J_p = [y_{p-1},y_p), p=1,\dots,m$ 
with a permutation $\sigma$ on $\{1,2,\dots,m \}$
such that
\begin{equation*}
f(x_{p-1}) = y_{\sigma(p)-1},
\qquad
f_{-}(x_p) = y_{\sigma(p-1)+1},
\qquad
 p=1,2,\dots,m
\end{equation*}
where
$
f_{-}(x_p) = \lim_{h \to 0+} f(x_p -h),
$
and $f$ is linear on $[x_{p-1}, x_p)$
with slope
$\frac{y_{\sigma(p)} -y_{\sigma(p)-1}}{x_p -x_{p-1}}$
for
$ p=1,2,\dots,m$.
Similarly to the preceding section, we have
\begin{prop}
A $\beta$-adic PL function for sofic $\beta$-shift naturally 
gives rise to $\beta$-adic pattern of squares for sofic $\beta$-shift.
\end{prop}
We may directly construct a $\beta$-adic PL function $f_T$ for sofic $\beta$-shift
from a $\beta$-adic table for sofic $\beta$-shift
\begin{math}
T=
\bigl[
\begin{smallmatrix}
{\nu(1)_{[i_1]}} & {\nu(2)_{[i_2]}} & \cdots & {\nu(m)_{[i_m]}} \\
{\mu(1)_{[i_1]}} & {\mu(2)_{[i_2]}} & \cdots & {\mu(m)_{[i_m]}} 
\end{smallmatrix}
\bigr]
\end{math}
as follows.
Put
$x_j = l(\nu(j+1)_{[i_j]}), 
\hat{y}_j = l(\mu(j+1)_{[i_j]}), 
j=0,1,\dots, m-1$.
Define $f_T$ by
$f_T(x_j) = \hat{y}_j, j=0, 1,\dots,m-1$
and
$f_T$ is linear on $[x_{j-1},x_j), j=1,2,\dots,m$
with slope
$\frac{r(\mu(j))-l(\mu(j))}{r(\nu(j))-l(\nu(j))} 
= \beta^{|\nu(j)|-|\mu(j)|}.$
The function $f_T$ yields a $\beta$-adic PL function
for sofic $\beta$-shift.

It is straightforward to see that 
the composition of two $\beta$-adic PL functions
for sofic $\beta$-shift is also a $\beta$-adic PL function for sofic $\beta$-shift.
Hence
the set of $\beta$-adic PL functions for sofic $\beta$-shift forms a group
under compositions.
We reach the following:
\begin{thm}\label{thm:soficPL}
The topological full group
$\Gamma_\beta$
for a sofic $\beta$-shift $(X_\beta,\sigma)$
is realized as
the group of all $\beta$-adic PL functions for sofic $\beta$-shift.
\end{thm}

\section{Classification of the topological ful groups $\Gamma_\beta$}
In this section we will classify the groups 
$\Gamma_\beta$
for SFT $\beta$-shifts and sofic $\beta$-shifts.
We will first classify $\Gamma_\beta$
for SFT $\beta$-shifts.

1. SFT case:
\begin{prop}
Suppose that
the $\beta$-shift
$(X_\beta,\sigma)$ is a shift of finite type 
such that the  $\beta$-expansion
of $1$ is
$
1 = 
\frac{\eta_1}{\beta} +\frac{\eta_2}{\beta^2} +\cdots + \frac{\eta_n}{\beta^n}. 
$
Set
\begin{align*}
T_i & = S_{i-1} \quad \text{ for } i=1,\dots, \eta_1, \\
T_{\eta_1 +i} & = S_{\eta_1}S_{i-1} \quad \text{ for } i=1,\dots, \eta_2, \\
T_{\eta_1+\eta_2 +i} 
& = S_{\eta_1}S_{\eta_2}S_{i-1} \quad \text{ for } i=1,\dots, \eta_3, \\
\vdots
&  \\
T_{\eta_1 + \eta_2 +\cdots + \eta_{n-1}+i} 
& = S_{\eta_1}S_{\eta_2} \cdots S_{\eta_{n-1}}S_{i-1} 
\quad \text{ for } i=1,\dots, \eta_n. 
\end{align*}
Define the $C^*$-subalgebras 
$\widehat{\mathcal{O}}_\beta,  \widehat{\mathcal{D}}_\beta$
 of $\OB$ by
\begin{align*}
\widehat{\mathcal{O}}_\beta 
& = C^*(T_i ; i=1,2,\dots, \eta_1 +\eta_2+\cdots+\eta_n),\\
\widehat{\mathcal{D}}_\beta 
& = C^*(T_\mu T_\mu^* ; 
\mu = (\mu_1,\dots,\mu_m), \mu_i = 1,2, \dots,\eta_1 +\eta_2+\cdots+\eta_n).
\end{align*}
Then 
the $C^*$-algebras 
$\widehat{\mathcal{O}}_\beta$
and
$\widehat{\mathcal{D}}_\beta$
coincide with
$\OB$ and ${\mathcal{D}}_\beta$
respectively,
and are isomorphic to
the Cuntz algebra 
${\mathcal{O}}_{\eta_1 +\eta_2+\cdots+\eta_n}$
and the canonical Cartan subalgebra
${\mathcal{D}}_{\eta_1 +\eta_2+\cdots+\eta_n}$
respectively, that is
\begin{equation*}
\widehat{\mathcal{O}}_\beta 
 = {\mathcal{O}}_\beta 
 = {\mathcal{O}}_{\eta_1 +\eta_2+\cdots+\eta_n},
 \qquad
\widehat{\mathcal{D}}_\beta 
 = {\mathcal{D}}_\beta 
 = {\mathcal{D}}_{\eta_1 +\eta_2+\cdots+\eta_n}.
\end{equation*}
\end{prop}
\begin{pf}
It is direct to see that the operators
$T_1,T_2, \dots, T_{\eta_1 +\eta_2+\cdots+\eta_n}$
are all isometries.
We then have
\begin{align*}
\sum_{i=1}^{\eta_1}T_i T_i^* 
& =  \sum_{j=0}^{\eta_1-1} S_j S_j^* 
= 1- S_{\eta_1}S_{\eta_1}^*, \\ 
\sum_{i=\eta_1+1}^{\eta_1+\eta_2}T_i T_i^* 
& =  \sum_{j=0}^{\eta_2-1} S_{\eta_1}S_j S_j^*S_{\eta_1}^* 
= S_{\eta_1}(1- S_{\eta_2}S_{\eta_2}^*)S_{\eta_1}^*, \\
&  \vdots \\
\sum_{i=\eta_1+\eta_2 + \cdots+\eta_{n-1}+1}^{\eta_1+\eta_2 + \cdots+\eta_n}
T_i T_i^* 
& =  \sum_{j=0}^{\eta_n-1} S_{\eta_1}S_{\eta_2} \cdots S_{\eta_{n-1}}
S_j S_j^*S_{\eta_{n-1}}^* S_{\eta_1}^*\cdots S_{\eta_2}^*S_{\eta_1}^* \\
& = 
S_{\eta_1}S_{\eta_2} \cdots S_{\eta_{n-1}}(1-
S_{\eta_n} S_{\eta_n}^*) S_{\eta_{n-1}}^* \cdots S_{\eta_2}^*S_{\eta_1}^*\\
& = 
S_{\eta_1}S_{\eta_2} \cdots S_{\eta_{n-1}} S_{\eta_{n-1}}^* \cdots S_{\eta_2}^*S_{\eta_1}^*. 
\end{align*}
It follows that
\begin{align*}
& \sum_{i=1}^{\eta_1+\eta_2 + \cdots+\eta_n}T_i T_i^* \\
=
& \sum_{i=1}^{\eta_1}T_i T_i^* 
+ \sum_{i=\eta_1+1}^{\eta_1+\eta_2}T_i T_i^*
+\cdots +
\sum_{i=\eta_1+\eta_2 + \cdots+\eta_{n-1}+1}^{\eta_1+\eta_2 + \cdots+\eta_n}
T_i T_i^* \\ 
= 
&1- S_{\eta_1}S_{\eta_1}^*
+ S_{\eta_1}(1- S_{\eta_2}S_{\eta_2}^*)S_{\eta_1}^*
+ \cdots +
S_{\eta_1}S_{\eta_2} \cdots S_{\eta_{n-1}} S_{\eta_{n-1}}^* \cdots S_{\eta_2}^*S_{\eta_1}^*  =1. 
\end{align*}
Hence the $C^*$-algebra 
$\widehat{\mathcal{O}}_\beta$
is isomorphic to the Cuntz algebra
${\mathcal{O}}_{\eta_1 + \eta_2+ \cdots + \eta_n}$.
The inclusion relation
$\widehat{\mathcal{O}}_\beta \subset \mathcal{O}_\beta$
is clear.
To show the converse 
inclusion relation
$\mathcal{O}_\beta \subset \widehat{\mathcal{O}}_\beta$,
it suffices to prove that 
the partial isometry
$S_{\eta_1}$
belongs to
the algebra
$\widehat{\mathcal{O}}_\beta$.
By the equality
\begin{equation*}
\varphi(S_{\eta_1}^*S_{\eta_1})
= \beta - \eta_1 
= 
\frac{\eta_2}{\beta} +\frac{\eta_3}{\beta^2} +\cdots + 
\frac{\eta_n}{\beta^{n-1}}, 
\end{equation*}
we have
\begin{align*}
S_{\eta_1}^*S_{\eta_1}
= & \sum_{j=0}^{\eta_2-1}S_j S_j^* + 
    \sum_{j=0}^{\eta_3-1} S_{\eta_2}S_j S_j^* S_{\eta_2}^* + \cdots \\
 +& \sum_{j=0}^{\eta_n-1}
    S_{\eta_2}S_{\eta_3}\cdots S_{\eta_{n-1}}S_j 
    S_j^* S_{\eta_n-1}^* \cdots S_{\eta_3}^*S_{\eta_2}^*  
\end{align*}
so that
\begin{align*}
S_{\eta_1} 
= & S_{\eta_1}S_{\eta_1}^*S_{\eta_1}\\
= & \sum_{j=0}^{\eta_2-1} T_{\eta_1+j+1}S_j^* + 
    \sum_{j=0}^{\eta_3-1} T_{\eta_1+\eta_2+j+1}(S_{\eta_2} S_j)^* + \cdots \\
 +& \sum_{j=0}^{\eta_n-1}
    T_{\eta_1+\eta_2+\cdots+\eta_{n-1}+j+1} 
    (S_{\eta_2}S_{\eta_3}\cdots S_{\eta_{n-1}}S_j)^*.  
\end{align*}
Denote by $\eta_0$ the empty word.
The following set $W_\beta$ of the words
\begin{equation*}
W_\beta =\{
(\eta_2, \eta_3, \dots, \eta_{m-1}, i)
\mid i=0,1,\dots,\eta_m-1, m=1,2,\dots,n \}
\end{equation*}
are all admissible words of $X_\beta$.
By cutting a word in the subwords begining with $\eta_1$,
one easily sees that any admissible word of $X_\beta$ is decomposed into a product of some of words of the following set   
\begin{equation*}
C_\beta =\{
(\eta_1, \eta_2, \dots, \eta_{m-1}, j) 
\mid j=0,1,\dots,\eta_m-1, m=1,2,\dots,n \}.
\end{equation*}
Hence any word of $W_\beta$ is a product of some of words of $C_\beta$.
This implies that the operators
\begin{equation*}
S_{\eta_2}S_{\eta_3}\cdots S_{\eta_{m-1}}S_j,
\quad j=0,1,\dots,\eta_m-1, m=1,2,\dots,n
\end{equation*}
are products of some of 
$
T_1, T_2, \dots, T_{\eta_1 + \eta_2 + \cdots + \eta_n}.
$
Therefore $S_{\eta_1}$ is written as a product of 
$T_i, T_i^*, i=1,2, \dots,\eta_1 + \eta_2 + \cdots + \eta_n.$
This shows that
$\mathcal{O}_\beta \subset \widehat{\mathcal{O}}_\beta$.
The equality $\DB =\widehat{\mathcal{D}}_\beta$ 
is direct.
\end{pf}
The above proposition implies that 
the SFT $\beta$-shift
$(X_\beta,\sigma)$ is continuously orbit equivalent to
the full $(\eta_1 +\eta_2+\cdots+\eta_n)$-shift
$(X_{\eta_1 +\eta_2+\cdots+\eta_n},\sigma)$
(\cite{MaPacific}, \cite{MatuiPLMS}, \cite{ReIMSB}).
Therefore we have
\begin{thm}\label{thm:HTfinite}
If the $\beta$-expansion of $1$ 
is finite such that 
\begin{equation*}
1 = 
\frac{\eta_1}{\beta} +\frac{\eta_2}{\beta^2} +\cdots + \frac{\eta_n}{\beta^n}, 
\end{equation*}
then the group
$\Gamma_\beta$ is 
isomorphic to the Higman-Thompson group
$V_{\eta_1 +\eta_2+\cdots+\eta_n}$.
\end{thm}
\begin{cor}
Let
$(X_\beta,\sigma)$ 
and
$(X_{\beta'},\sigma)$ 
be SFT $\beta$-shifts such that their finite $\beta$-expansions
of $1$ are
\begin{equation*}
1 = 
\frac{\eta_1}{\beta} +\frac{\eta_2}{\beta^2} +\cdots + \frac{\eta_n}{\beta^n} 
= 
\frac{\eta'_1}{\beta'} +\frac{\eta'_2}{{\beta'}^2} +\cdots 
+ \frac{\eta'_{n'}}{{\beta'}^{n'}}
\end{equation*}
respectively.
Then the followings are equivalent:
\begin{enumerate}
\renewcommand{\labelenumi}{(\roman{enumi})}
\item The groups
$\Gamma_\beta$ and 
$\Gamma_{\beta'}$ are isomorphic.
\item
The Cuntz-algebras
${\mathcal O}_{\eta_1 + \eta_2 +\cdots + \eta_n}$
and
${\mathcal O}_{\eta'_1 + \eta'_2 +\cdots + \eta'_{n'}}$
are isomorphic.
\item
$\eta_1 + \eta_2 +\cdots + \eta_n
= \eta'_1 + \eta'_2 +\cdots + \eta'_{n'}$.
\end{enumerate}
\end{cor}
\begin{pf}
The implication (iii) $\Longrightarrow$ (ii) is trivial,
and its converse (ii) $\Longrightarrow$ (iii) is well-known
(\cite{CuntzCMP}, \cite{CuntzAnnMath}).
Assume that
the groups
$\Gamma_\beta$ and 
$\Gamma_{\beta'}$ are isomorphic.
By \cite{MaPre2012} or more generally \cite{MatuiPre2012},
the $C^*$-algebras
$C^*_r(G_\beta)$
and
$C^*_r(G_{\beta'})$
of the groupoids 
$G_\beta$ and $G_{\beta'}$
associated to their respective  shifts
$(X_\beta,\sigma)$
and
$(X_{\beta'},\sigma)$
of finite type are isomorphic.
Since  
$C^*_r(G_\beta)=\OB$
and
$C^*_r(G_{\beta'})={\mathcal{O}}_{\beta'}$,
Proposition 7.1 implies (ii), so that
the implication (i) $\Longrightarrow$ (ii) holds. 
The implication
(iii) $\Longrightarrow$ (i) 
is a direct consequence of the above theorem.
\end{pf}

\medskip

2. Sofic case:

Assume that
the $\beta$-shift $X_\beta$ is sofic.
Put
\begin{equation*}
k_\beta = \min\{ k \in {\mathbb{N}} \mid \A_k = \A_{k+1}\},
\qquad
K_\beta = k_\beta +1.
\end{equation*}
Hece
$\A_{k_\beta} = \A_{k_\beta+1} =\cdots = \AB$
and
$\dim \AB = K_\beta$.
There exists $l \in {\mathbb{N}}$ 
with $0 < l \le k_\beta$
such that 
\begin{equation}
a_{\xi_1 \cdots \xi_{K_\beta}} = a_{\xi_1 \cdots \xi_l}
\quad
\text{ and hence }
\quad 
d(1,\beta) = \xi_1 \cdots \xi_l \Dot{\xi}_{l+1} \cdots \Dot{\xi}_{K_\beta}. 
\label{eqn:kbl}
\end{equation}
Let
$E_1,\dots, E_{K_\beta}$ be the minimal projections of $\AB$
as in the preceding section so that
\begin{equation}
\AB = {\mathbb{C}}E_1 \oplus \cdots \oplus {\mathbb{C}}E_{K_\beta}. 
\label{eqn:AB}
\end{equation}
Define a labeled graph ${\mathcal{G}}_\beta$
over $\Sigma = \{0,1,\dots,N-1\}$
with vertex set 
$\{ v_1, v_2,\dots, v_{K_\beta}\}$
corresponding to the minimal projections
$E_1,\dots, E_{K_\beta}$
in the following way.
Define a labeled edge from $v_i$ to $v_j$ 
labeled $\alpha \in \Sigma$ if
$S_\alpha^* E_i S_\alpha \ge E_j$.
The underlying directed graph
of the labeled graph ${\mathcal{G}}_\beta$ is   
denoted by $G_\beta = ({\mathcal{V}}_\beta, {\mathcal{E}}_\beta)$
with labeling map 
$\lambda : {\mathcal{E}}_\beta \longrightarrow \Sigma$.
Let 
${\mathcal{M}}_\beta$ 
be the $K_\beta \times K_\beta$
symbolic matrix of 
${\mathcal{G}}_\beta$
and 
$M_\beta$  the $K_\beta \times K_\beta$ 
nonnegative matrix obtained from
${\mathcal{M}}_\beta$ by putting all the symbols equal to $1$.
For an edge $e \in {\mathcal{E}}_\beta$,
put a partial isometry
$s_e = S_{\lambda(e)}E_{t(e)}$
in the $C^*$-algebra $\OB$,
where
$\lambda(e) \in \Sigma$ 
denotes the letter of the label of $e$
and 
$t(e) \in \{1,2,\dots, K_\beta\}$ 
denotes the number of the terminal vertex $v_{t(e)}$  of $e$.
Define the 
$|{\mathcal{E}}_\beta|\times |{\mathcal{E}}_\beta|$ matrix
$B_\beta = [B_\beta(e,f)]_{e,f \in {\mathcal{E}}_\beta}$ 
with entries in $\{0,1\}$ by
\begin{equation*}
B_\beta(e,f)=
\begin{cases}
1 & \text{ if } t(e) = s(f),\\
0 & \text{ if } t(e) \neq  s(f).
\end{cases}
\end{equation*}
We have the following lemma.
\begin{lem}
The partial isometries 
$s_e, e \in {\mathcal{E}}_\beta$
satisfy the following relations:
\begin{equation*}
\sum_{e \in {\mathcal{E}}_\beta} s_e s_e^* =1,
\qquad
s_e^* s_e = \sum_{f \in {\mathcal{E}}_\beta} B_\beta(e,f) s_f s_f^*.
\end{equation*}
Hence the $C^*$-algebra 
$C^*(s_e; e\in  {\mathcal{E}}_\beta)$
generated by $s_e, e \in {\mathcal{E}}_\beta$
is isomorphic to the Cuntz-Krieger algebeta
${\mathcal{O}}_{B_\beta}$.
\end{lem}
\begin{pf}
We see the identities 
\begin{equation*}
1 = \sum_{i=1}^{K_\beta} E_i
  = \sum_{i=1}^{K_\beta} \sum_{\alpha =0}^{N-1}
     S_\alpha S_\alpha^* E_i S_\alpha S_\alpha^*.
\end{equation*}
For an edge $e \in {\mathcal{E}}_\beta$,
denote by 
$s(e) \in \{1,2,\dots,K_\beta\}$
the number of the source vertex $v_{s(e)}$  of $e$.
 The projection 
$S_\alpha^* E_i S_\alpha$ is not zero if and only if 
there exists 
$e \in {\mathcal{E}}_\beta$
such that 
$\alpha = \lambda(e)$ and 
$i = s(e)$.
Hence we have
\begin{equation*}
S_\alpha^* E_i S_\alpha
=
\sum
\begin{Sb}
e \in {\mathcal{E}}_\beta,\\
\alpha = \lambda(e),i = s(e)
\end{Sb}
E_{t(e)}
\end{equation*}
so that
\begin{equation*}
1 = \sum_{i=1}^{K_\beta} \sum_{\alpha =0}^{N-1}
    \sum
\begin{Sb}
e \in {\mathcal{E}}_\beta,\\
\alpha = \lambda(e),i = s(e)
\end{Sb}
S_\alpha E_{t(e)}S_\alpha^* 
 = \sum_{e \in {\mathcal{E}}_\beta} s_e s_e^*.
\end{equation*}
For an edge
$e \in {\mathcal{E}}_\beta$,
we see that
\begin{align*}
s_e^* s_e 
& = E_{t(e)} S_{\lambda(e)}^* S_{\lambda(e)} E_{t(e)} = E_{t(e)} \\
& = \sum_{\alpha =0}^{N-1} S_\alpha S_\alpha^* E_{t(e)} S_\alpha S_\alpha^* \\
&  = \sum_{\alpha =0}^{N-1} S_\alpha \cdot
    \sum
\begin{Sb}
f \in {\mathcal{E}}_\beta,\\
\alpha = \lambda(f),t(e) = s(f)
\end{Sb}
 E_{t(f)} \cdot S_\alpha^* 
 = \sum_{f \in {\mathcal{E}}_\beta} B_\beta(e,f) s_f s_f^*.
\end{align*}
\end{pf}

Denote by 
${\mathcal{D}}_{B_\beta}$
 the canonical Cartan subalgebra of ${\mathcal{O}}_{B_\beta}$
which is a $C^*$-subalgebra of ${\mathcal{O}}_{B_\beta}$
generated by projections
$s_{e_1}\cdots s_{e_n} s_{e_n}^*\cdots s_{e_1}^*$,
$ 
e_1,\dots,e_n \in  {\mathcal{E}}_\beta.
$
\begin{lem}
We have
\begin{equation*}
\OB = {\mathcal{O}}_{B_\beta},\qquad
\DB = {\mathcal{D}}_{B_\beta}.
\end{equation*}
\end{lem}
\begin{pf}
Since $s_e = S_{\lambda(e)} E_{t(e)}, e \in {\mathcal{E}}_\beta$,
we have 
$s_e \in \OB$
so that
the inclusion
${\mathcal{O}}_{B_\beta} \subset \OB$ is obvious.
For $\alpha \in \Sigma = \{0,1,\dots,N-1\}$,
$i=1,\dots,K_\beta$, 
we know that 
$S_\alpha E_i \ne 0$ 
if and only if 
$S_\alpha^* S_\alpha \ge E_i$,
which is equivalent 
to the condition that 
there exists an edge
$e \in {\mathcal{E}}_\beta$ such that
$\alpha = \lambda(e), i =t(e)$.
For $i=1,\dots,K_\beta$, 
take $e \in {\mathcal{E}}_\beta$ 
such that
$\alpha = \lambda(e), i =t(e)$.
We then have 
$
s_e^* s_e 
= E_{t(e)} = E_i
$ 
so that
$E_i \in {\mathcal{O}}_{B_\beta}$.
For $\alpha \in \Sigma$, we have
\begin{equation*}
 S_\alpha 
 = \sum_{i=1}^{K_\beta} S_\alpha E_i
 = \sum_{e \in {\mathcal{E}}_\beta, \alpha=\lambda(e)}
 S_{\lambda(e)} E_{t(e)}
= \sum_{e \in {\mathcal{E}}_\beta, \alpha=\lambda(e)}
 s_e
\end{equation*} 
so that
$S_\alpha \in {\mathcal{O}}_{B_\beta}$.
Therefore we have
the inclusion
$\OB \subset {\mathcal{O}}_{B_\beta}$ 
and hence
$\OB = {\mathcal{O}}_{B_\beta}$.

We will next show that 
$\DB = {\mathcal{D}}_{B_\beta}$.
We have
$s_e s_e^* = S_{\lambda(e)} E_{t(e)} S_{\lambda(e)}^* \in \DB$.
Suppose that 
$s_{e_1}\cdots s_{e_n} s_{e_n}^*\cdots s_{e_1}^* \in \DB$.
By the equality
$$
s_{e_0}s_{e_1}\cdots s_{e_n} s_{e_n}^*\cdots s_{e_1}^*s_{e_0}^*
=S_{\lambda(e_0)} E_{t(e_0)}s_{e_1}\cdots s_{e_n} s_{e_n}^*\cdots s_{e_1}^*
E_{t(e_0)}S_{\lambda(e_0)}^*
 \in \DB,
$$
the inclusion relation
${\mathcal{D}}_{B_\beta} \subset \DB$ holds
by induction.
Conversely suppose that
$S_\alpha E_i S_\alpha^*$ is not zero.
Take 
$e \in {\mathcal{E}}_\beta$
such that 
$ \alpha=\lambda(e), i = t(e)$
so that
$$
S_\alpha E_i S_\alpha^*
=
S_{\lambda(e)} E_{t(e)}S_{\lambda(e)}^*
= s_e s_e^*
$$
that belongs to
${\mathcal{D}}_{B_\beta}$.
Suppose next that 
$
S_{\mu_1 \cdots \mu_n}E_i S_{\mu_1 \cdots \mu_n}^*$ 
belongs to
$
{\mathcal{D}}_{B_\beta}
$
and
$
S_{\mu_0}S_{\mu_1 \cdots \mu_n}E_i S_{\mu_1 \cdots \mu_n}^* S_{\mu_0}^*
$
is not zero.
The labeled graph
${\mathcal{G}}_\beta$ is left-resolving
which means that
there uniquely exists
a finite sequence of edges
$e_1, e_2, \dots,e_n  \in {\mathcal{E}}_\beta$
for the vertex $v_i$
such that 
\begin{align*}
\lambda(e_p) & = \mu_p, 
\,
t(e_p) = s(e_{p+1})
\quad
\text{ for }
  p=1,\dots, n-1,\\
\lambda(e_n) & = \mu_n,  
\,
t(e_n) = i.
\end{align*} 
Put
$j = s(e_1)$
so that we have
$$
E_j \ge S_{\mu_1 \cdots \mu_n}E_i S_{\mu_1 \cdots \mu_n}^*,
\quad
E_j S_{\mu_1 \cdots \mu_n}E_i S_{\mu_1 \cdots \mu_n}^*E_j
= S_{\mu_1 \cdots \mu_n}E_i S_{\mu_1 \cdots \mu_n}^*.
$$
Take a unique edge
$e_0  \in {\mathcal{E}}_\beta$
such that
$\lambda(e_0) = \mu_0, t(e_0) = j$.
Hence 
$S_{\mu_0} E_j = s_{e_0}$.
It then follows that
\begin{align*}
S_{\mu_0}S_{\mu_1 \cdots \mu_n}E_i S_{\mu_1 \cdots \mu_n}^* S_{\mu_0}^*
& = 
S_{\mu_0}E_j S_{\mu_1\cdots \mu_n}E_i S_{\mu_1\cdots \mu_n}^*E_j S_{\mu_0}^*\\
& = 
s_{e_0}S_{\mu_1 \cdots \mu_n}E_i S_{\mu_1 \cdots \mu_n}^*s_{e_0}^*.
\end{align*}
As
$
S_{\mu_1 \cdots \mu_n}E_i S_{\mu_1 \cdots \mu_n}^* \in 
{\mathcal{D}}_{B_\beta},
$
we have
$
s_{e_0}S_{\mu_1 \cdots \mu_n}E_i S_{\mu_1 \cdots \mu_n}^*s_{e_0}^*
\in 
{\mathcal{D}}_{B_\beta}.
$
Thus 
the element
$
S_{\mu_0}S_{\mu_1 \cdots \mu_n}E_i S_{\mu_1 \cdots \mu_n}^* S_{\mu_0}^*
$
belongs to
$
{\mathcal{D}}_{B_\beta}.
$
By induction, we have    
$\DB \subset {\mathcal{D}}_{B_\beta}$
and hence
$\DB = {\mathcal{D}}_{B_\beta}$.
\end{pf}
A nonnegative square matrix $B$ is said to be strong shift equivalent to
a nonnegative square matrix $M$ 
if there exist nonnegative rectangular matrices $R$ and $S$ such that
$B=RS$ and $M =SR$ (cf. \cite{LM}). 
\begin{lem}
The matrix $B_\beta$  
is strong shift equivalent to 
the matrix $M_\beta$.
Hence we have
\begin{equation*}
\det(1-B_\beta) =\det(1-M_\beta).
\end{equation*}
\end{lem}
\begin{pf}
Note that 
$\dim {\mathcal{A}}_\beta = | \mathcal{V}_\beta | = K_\beta$.
Define 
a $|{\mathcal{E}}_\beta | \times |{\mathcal{V}}_\beta |$ matrix $R_\beta$
and  
a $|{\mathcal{V}}_\beta | \times |{\mathcal{E}}_\beta |$ matrix $S_\beta$
as follows.
\begin{equation*}
R_\beta(e,i) = 
\begin{cases}
1 & \text{ if } t(e) = v_i,\\
0 & \text{ otherwise },
\end{cases}
\qquad 
S_\beta(j, f) = 
\begin{cases}
1 & \text{ if } s(f) = v_j,\\
0 & \text{ otherwise }
\end{cases}
\end{equation*}
for
$e,f \in {\mathcal{E}}_\beta,
v_i, v_j \in {\mathcal{V}}_\beta
$
and
$i,j =1,\dots,K_\beta$.
It is direct to see that
\begin{equation*}
B_\beta = R_\beta S_\beta, \qquad
M_\beta = S_\beta R_\beta
\end{equation*}
and 
$\det(1-B_\beta) =\det(1-M_\beta).$
\end{pf}
Recall that  
$\varphi$ stands for the unique KMS state 
on the $C^*$-algebra $\OB$ 
under the gauge action.
It satisfies the identities
\begin{equation*}
\varphi(a_{\xi_1 \cdots \xi_j}) =
      \beta^j- \xi_1 \beta^{j-1} - \xi_2 \beta^{j-2} - 
\cdots -\xi_{j-1}\beta - \xi_j, \quad  j =1,\dots, K_{\beta}. 
\end{equation*}
By \eqref{eqn:AB}, 
the $K_0$-group $K_0(\A_{k_{\beta}})$
of the algebra $\A_{k_{\beta}}$
is generated by the classes of the 
minimal projections 
$E_1,\dots, E_{K_\beta}$ of 
$\A_{k_{\beta}}(=\A_{\beta})$
so that $K_0(\A_{k_{\beta}})$
is isomorphic to
${\mathbb{Z}}^{K_\beta}$.
Since a minimal projection $E_i$
is of the form
$
a_{\xi_1 \cdots \xi_{p_i}} -
a_{\xi_1 \cdots \xi_{q_i}},
$
the following correspondence
\begin{align*}
[1] \in K_0(\A_{k_{\beta}}) & \longrightarrow (1,0,0,\dots,0) 
\in {\mathbb{Z}} \oplus \beta{\mathbb{Z}} \oplus \cdots \oplus \beta^{k_\beta}{\mathbb{Z}} ,\\ 
[a_{\xi_1}] \in K_0(\A_{k_{\beta}}) & \longrightarrow (-\xi_1, 1, 0,\dots,0) 
\in {\mathbb{Z}} \oplus \beta{\mathbb{Z}} \oplus \cdots \oplus \beta^{k_\beta}{\mathbb{Z}} ,\\
[a_{\xi_1\cdots \xi_j}] \in K_0(\A_{k_{\beta}}) & \longrightarrow 
(-\xi_j,-\xi_{j-1}, \dots, -\xi_2, -\xi_1, 1, 0,\dots,0) 
\in {\mathbb{Z}} \oplus \beta{\mathbb{Z}} \oplus \cdots \oplus \beta^{k_\beta}{\mathbb{Z}} 
\end{align*}
for $j =1,\dots, K_{\beta}$ 
yields an isomorphism from
$K_0(\A_{k_{\beta}})$ to
${\mathbb{Z}} \oplus \beta{\mathbb{Z}} \oplus \cdots \oplus \beta^{k_\beta}{\mathbb{Z}} $
as a group, which we denote by $\varPhi$. 
By \eqref{eqn:kbl}, we have
\begin{align*}
  & \beta^{k_{\beta}+1}- \xi_1 \beta^{k_{\beta}} - \xi_2 \beta^{k_{\beta}-1} - 
\cdots 
   -\xi_{k_{\beta}}\beta - \xi_{k_{\beta}+1} \\
   = 
  & \beta^l- \xi_1 \beta^{l-1} - \xi_2 \beta^{l-2} - 
\cdots -\xi_{l-1}\beta - \xi_l
\end{align*}
so that $\beta$ is a solution of a monic polynomial of degree $K_{\beta}$.
We denote this polynomial by
\begin{equation*}
   \beta^{k_{\beta} +1}- \eta_1 \beta^{k_{\beta}} - \eta_2 \beta^{k_{\beta}-1} - 
\cdots 
   -\eta_{k_{\beta}}\beta - \eta_{k_{\beta}+1} 
   = 0.  
\end{equation*}
Then we have 
\begin{equation}
\eta_1 +\eta_2+ \cdots + \eta_{k_\beta} +\eta_{k_\beta +1} 
= \xi_{l+1} + \xi_{l+2} + \cdots + \xi_{k_\beta +1} +1. \label{eqn:etaxi}
\end{equation} 
\begin{lem}[{\cite[Lemma 4.8]{KMW}}]
The following diagram is commutative:
\begin{equation*}
\begin{CD}
{\mathbb{Z}}^{K_\beta}  @>{M_\beta} >> {\mathbb{Z}}^{K_\beta} \\
@|                     @|                                  \\
K_0(\A_{k_{\beta}})  @>{\lambda_{\beta *}} >> K_0(\A_{k_{\beta}})  \\
@V{\varPhi}VV              @V{\varPhi}VV           \\
{\mathbb{Z}} \oplus \beta {\mathbb{Z}} \oplus \cdots 
\oplus \beta^{k_{\beta}}{\mathbb{Z}} 
@>{\tau} >>
{\mathbb{Z}} \oplus \beta{\mathbb{Z}} \oplus \cdots 
\oplus \beta^{k_{\beta}}{\mathbb{Z}} 
\end{CD}
\end{equation*}
where 
$
\lambda_{\beta *}
$
is the endomorphism of $K_0(\A_{\beta})$ 
induced from the map
$\lambda_{\beta} : \A_{k_{\beta}} \rightarrow \A_{{k_{\beta}}+1} (= \A_{\beta})$
defined by
\begin{equation*}
\lambda_\beta (a) = \sum_{\alpha=0}^{N-1} S_\alpha^* a S_\alpha
\qquad \text{ for }
a \in {\mathcal{A}}_\beta
\end{equation*}
and $\tau$ is an endomorphism of
${\mathbb{Z}} \oplus \beta{\mathbb{Z}} \oplus \cdots 
\oplus \beta^{k_{\beta}}{\mathbb{Z}}
$
defined by
\begin{align*}
\tau(m_0,m_1, \dots,m_{k_{\beta}-1}, 0) 
& = (0,m_0,m_1, \dots,m_{k_{\beta}-1}), 
\quad m_i \in \Bbb{Z},\\
\tau(0, \dots ,0,1) 
& = (\eta_{k_{\beta}+1},\eta_{k_{\beta}}, \dots, \eta_2, \eta_1).
\end{align*}
\end{lem}
Set the
$(k_{\beta}+1) \times (k_{\beta}+1)$ matrix 
$$
L_{\beta} = 
\begin{bmatrix}
              &          &             &  \eta_{{k_{\beta}} +1} \\
      1       &          &             &  \eta_{k_{\beta}}      \\
              &  \ddots  &             &  \vdots          \\
              &          &     1       &  \eta_1          
\end{bmatrix}
$$
where the blanks denote zeros.
The matrix $L_{\beta}$ acts from the left hand side of the transpose
$(m_0,m_1,\dots,m_{k_\beta})^t$
of
$(m_0,m_1,\dots,m_{k_\beta})$
so that it represents the homomorphism $\tau$.
The characteristic polynomial of $L_{\beta}$ is 
\begin{equation*}
  \det(t -L_{\beta}) = 
     t^{{k_{\beta}} +1} - \eta_1 t^{k_{\beta}} - \eta_2 t^{k_{\beta} -1} -
  \cdots - \eta_{k_\beta}t - \eta_{k_\beta +1}
\end{equation*}
and
the number $\beta$ is one of the eigenvalues of the transpose of  
$L_{\beta}$ with eigenvector
$[1,\beta, \beta^2, \dots, \beta^{k_\beta}].$
Hence we have
\begin{cor}
$\det(1-B_\beta)  
= \det(1-L_\beta) 
= 1-\eta_1-\eta_2-\cdots - \eta_{k_\beta} - \eta_{k_\beta +1} <0.
$
\end{cor}
\begin{prop}
There exists an isomorphism
$\varPhi$ from the Cuntz-Krieger algebra
${\mathcal{O}}_{B_\beta}$ onto the 
 the Cuntz algebra 
$
{\mathcal{O}}_{\xi_1 + \cdots +\xi_{k_\beta +1} +1}
$
such that
$\varPhi({\mathcal{D}}_{B_\beta}) 
={\mathcal{D}}_{\xi_1 + \cdots +\xi_{k_\beta +1} +1}.
$
Therefore their topological  full groups
$\Gamma_{B_\beta}$
and
$
\Gamma_{\xi_1 + \cdots +\xi_{k_\beta +1} +1}
$
are isomorphic.
\end{prop}
\begin{pf}
We have already known that
$\OB$ is isomorphic to
${\mathcal{O}}_{\xi_1 + \cdots +\xi_{k_\beta +1} +1}
$
by \cite{KMW}.
By the preceding lemma,
we know that
$\OB ={\mathcal{O}}_{B_\beta}$
and
$\DB ={\mathcal{D}}_{B_\beta}$
so that
${\mathcal{O}}_{B_\beta}$
is isomorphic to
${\mathcal{O}}_{\xi_1 + \cdots +\xi_{k_\beta +1} +1}.
$
By the preceding lemma with \eqref{eqn:etaxi},
we see that
\begin{align*}
\det(1-B_\beta)
& = 1-\eta_1-\eta_2-\cdots - \eta_{k_\beta} - \eta_{k_\beta +1} \\
& = 1- (\xi_{l+1} + \cdots +\xi_{k_\beta +1} +1).
\end{align*} 
Hence the topological Markov shift
$(X_{B_\beta}, \sigma)$
is continuously orbit equivalent to
the full shift
$(X_{\xi_{l+1} + \cdots +\xi_{k_\beta +1} +1}, \sigma)$
by \cite{MaPAMS}.
Thus their topological full groups 
$\Gamma_{B_\beta}$
and
$
\Gamma_{\xi_{l+1} + \cdots +\xi_{k_\beta +1} +1}
$
are isomorphic.
\end{pf}
\begin{thm}\label{thm:HTsofic}
Suppose that 
$(X_\beta,\sigma)$ is sofic
such that the $\beta$-expansion of $1$ is
\begin{equation*}
d(1,\beta) = \xi_1  \cdots \xi_l \overset{\cdot}{\xi}_{l+1}  \cdots \overset{\cdot}{\xi}_{k +1}.
\end{equation*}
Then 
there exists an isomorphism
$\Phi$ 
from
$\OB$ to
$
 {\mathcal{O}}_{\xi_{l+1} + \cdots +\xi_{k +1} +1}
$
such that
$\Phi(\DB)
={\mathcal{D}}_{\xi_{l+1} + \cdots +\xi_{k +1} +1}.
$
Therefore their topological full groups
$
\Gamma_{\beta}
$
and
$
 \Gamma_{\xi_{l+1} + \cdots +\xi_{k +1} +1}
$
are isomorphic.
This implies that
the group
$\Gamma_\beta$  
is isomorphic to the Higman-Thompson group
$
V_{\xi_{l+1} + \cdots +\xi_{k +1} +1}.
$
\end{thm}

\medskip

3. Non-sofic case:

\begin{thm}\label{thm:HTnonsofic}
If $1 < \beta \in {\mathbb{R}}$
is not ultimately periodic,
then the group
$\Gamma_\beta$ is not isomorphic to any of 
the Higman-Thompson group $V_n, 1 < n \in {\mathbb{N}}$. 
\end{thm}
\begin{pf}
By Proposition \ref{prop:pureinf},
the groupoid $G_\beta$ is an essentially principal, purely infinite, minimal groupoid.
Suppose that
$\Gamma_\beta$ is isomorphic to
one of the Higman-Thompson group $V_n$ for some $n \in {\mathbb{N}}$.
Since
$V_n$ is isomorphic to the topological full group 
$\Gamma_n$ of the groupoid $G_n$ for the full $n$-shift. 
By H. Matui \cite{MatuiPre2012}, 
the groupoid
$G_\beta$ is isomorphic to
$G_n$.
By J. Renault \cite[Theorem 4.11]{Renault},
there exists an isomorphism
$\varPhi$ from 
$C^*_r(G_\beta)$ to $C^*_r(G_n)$.
The $C^*$-algebra
$C^*_r(G_\beta)$ is isomorphic to $\OB$,
and
 $C^*_r(G_n)$ is isomorphic to the Cuntz algebra ${\mathcal{O}}_n$.
Since $\beta$ is not ultimately periodic,
we know that $K_0(\OB) ={\mathbb{Z}}$
by \cite[Theorem 4.12]{KMW},
which is
a contradiction to the fact
$K_0({\mathcal{O}}_n) = {\mathbb{Z}}/(1-n){\mathbb{Z}}$.   
\end{pf}

\end{document}